\documentclass[brochure,english,11pt]{bourbaki}
\usepackage[matrix,arrow]{xy}
\usepackage{amssymb,amsfonts,amsmath,footnote}
\usepackage[english]{babel}
\date{Juin 2011}
\bbkannee{63\`eme ann\'ee, 2010-2011}
\bbknumero{1039}
\title{Invariant percolation and measured theory of nonamenable groups}
\subtitle{after Gaboriau-Lyons, Ioana, Epstein}
\author{Cyril HOUDAYER}
\address{\'Ecole Normale Sup\'erieure de Lyon\\
U.M.P.A.\\
CNRS UMR 5669\\
46 all\'ee d'Italie\\
F--69364 Lyon Cedex 07}
\email{cyril.houdayer@ens-lyon.fr}

\thanks{Research partially supported by ANR grant AGORA}

\newtheorem*{theorem}{Theorem}
\newtheorem*{corollary}{Corollary}
\newtheorem*{claim}{Claim}

\newtheorem{newstep}{Step}

\newcommand{\R}{\mathbf{R}}
\newcommand{\C}{\mathbf{C}}
\newcommand{\Z}{\mathbf{Z}}
\newcommand{\F}{\mathbf{F}}

\newcommand{\N}{\mathbf{N}}

\newcommand{\Id}{\operatorname{Id}}
\newcommand{\Ad}{\operatorname{Ad}}

\newcommand{\Tr}{\operatorname{Tr}}

\newcommand{\Aut}{\operatorname{Aut}}
\newcommand{\SL}{\operatorname{SL}}

\newcommand{\range}{\operatorname{range}}

\newcommand{\op}{\operatorname{op}}

\newcommand{\fin}{\operatorname{fin}}

\newcommand{\Stab}{\operatorname{Stab}}

\newcommand{\weak}{\operatorname{weak}}
\newcommand{\Cay}{\operatorname{Cay}}
\newcommand{\finite}{\operatorname{finite}}
\newcommand{\graph}{\operatorname{graph}}

\newcommand{\ME}{\operatorname{ME}}
\newcommand{\forest}{\operatorname{Forest}}
\newcommand{\cl}{\operatorname{cl}}
\newcommand{\QN}{\operatorname{QN}}
\newcommand{\coInd}{\operatorname{coInd}}
\newcommand{\cost}{\operatorname{cost}}
\newcommand{\Leb}{\operatorname{Leb}}

\begin{document}

\maketitle

\section{Introduction}

The notion of amenability was introduced in 1929 by J. von Neumann \cite{vonNeumann} in order to explain the Banach-Tarski paradox. A countable discrete group $\Gamma$ is {\em amenable} if there exists a left-invariant mean $\varphi : \ell^\infty(\Gamma) \to \C$. The class of amenable groups is stable under subgroups, direct limits, quotients and the free group $\F_2$ on two generators is not amenable. Knowing whether or not the class of amenable groups coincides with the class of groups without a nonabelian free subgroup became known as von Neumann's problem. It was solved in the negative by Ol'shanskii \cite{olshanski}. Adyan \cite{adyan} proved that the free Burnside groups $B(m, n)$ with $m$ generators, of exponent $n$ ($n \geq 665$ and odd) are nonamenable. Ol'shanskii and Sapir \cite{olshanski-sapir} also constructed examples of finitely presented nonamenable groups without a nonabelian free subgroup.

Two free ergodic probability measure-preserving (pmp) actions $\Gamma \curvearrowright (X, \mu)$ and $\Lambda \curvearrowright (Y, \nu)$ of countable discrete groups on nonatomic standard probability spaces are {\em orbit equivalent} (OE) if they induce the same orbit equivalence relation, that is, if there exists a pmp Borel isomorphism $\Delta : (X, \mu) \to (Y, \nu)$ such that $\Delta(\Gamma x) = \Lambda \Delta(x)$, for $\mu$-almost every $x \in X$. Despite the fact that the group $\Z$ admits uncountably many non-conjugate free ergodic pmp actions, Dye \cite{{dye1}, {dye2}} proved the surprising result that any two free ergodic pmp actions of $\Z$ are orbit equivalent. Moreover, Ornstein and Weiss \cite{ow}  (see also \cite{cfw}) proved that any free ergodic pmp action $\Gamma \curvearrowright (X, \mu)$ of any infinite amenable group is always orbit equivalent to a free ergodic pmp $\Z$-action on $(X, \mu)$. On the other hand, results of \cite{{schmidt}, {cw}, {hjorth-T}} imply that any nonamenable group has at least two non-OE free ergodic pmp actions. These results lead to a satisfying characterization of amenability: an infinite countable discrete group $\Gamma$ is amenable if and only if $\Gamma$ admits exactly one free ergodic pmp action up to OE.

\subsection*{Measurable-group-theoretic solution to von Neumann's problem} The first result we discuss in this paper is a positive answer to von Neumann's problem in the framework of measured group theory, due to Gaboriau and Lyons \cite{gaboriau-lyons}. Measured group theory is the study of countable discrete groups $\Gamma$ through their pmp actions $\Gamma \curvearrowright (X, \mu)$. We refer to \cite{gaboriau-icm} for a recent survey on this topic.

To any free pmp action $\Gamma \curvearrowright (X, \mu)$, one can associate the {\em orbit equivalence relation} $\mathcal R(\Gamma \curvearrowright X) \subset X \times X$ defined by
$$(x, y) \in \mathcal R(\Gamma \curvearrowright X) \Longleftrightarrow \exists g \in \Gamma, y = g x.$$
For countable discrete groups $\Gamma$ and $\Lambda$, we say that $\Lambda$ is a {\em measurable subgroup} of $\Gamma$ and denote $\Lambda <_{\ME} \Gamma$ if there exist two free ergodic pmp actions $\Gamma \curvearrowright (X, \mu)$ and $\Lambda \curvearrowright (X, \mu)$ such that $\mathcal R(\Lambda \curvearrowright X) \subset \mathcal R(\Gamma \curvearrowright X)$.
Denote by $\Leb$ the Lebesgue measure on the interval $[0, 1]$ and let $\Gamma \curvearrowright ([0, 1], \Leb)^\Gamma$ be the Bernoulli shift. Gaboriau and Lyons \cite{gaboriau-lyons} obtained the following remarkable result.
\begin{theorem}
Let $\Gamma$ be any nonamenable countable discrete group. Then there exists a free ergodic pmp action $\F_2 \curvearrowright ([0, 1], \Leb)^\Gamma$ such that 
$$\mathcal R(\F_2 \curvearrowright [0, 1]^\Gamma) \subset \mathcal R(\Gamma \curvearrowright [0, 1]^\Gamma).$$
\end{theorem}

In particular, we get that $\F_2 <_{\ME} \Gamma$. This theorem has important consequences in the theory of group von Neumann algebras.

\begin{corollary}
Let $\Gamma, H$ be countable discrete groups such that $\Gamma$ is nonamenable and $H$ is infinite. Then the von Neumann algebra $L(H \wr \Gamma)$ of the wreath product group $H \wr \Gamma := (\bigoplus_{\Gamma} H) \rtimes \Gamma$ contains a copy of the von Neumann algebra $L(\F_2)$ of the free group. 
\end{corollary}

The proof of Gaboriau and Lyons' result goes in two steps that we explain below. We refer to Section \ref{relations} for background material on pmp equivalence relations.

The {\bf first step} consists in finding a subequivalence relation $\mathcal R \subset \mathcal R(\Gamma \curvearrowright [0, 1]^\Gamma)$ such that $\mathcal R$ is ergodic treeable and non-hyperfinite. This is a difficult problem in general. By Zimmer's result \cite[Proposition 9.3.2]{zimmer}, it is known that $\mathcal R(\Gamma \curvearrowright [0, 1]^\Gamma)$ contains an ergodic hyperfinite subequivalence relation. When $\Gamma$ is finitely generated, another way to obtain subequivalence relations of $\mathcal R(\Gamma \curvearrowright [0, 1]^\Gamma)$ is by considering invariant {\em percolation} processes on the Cayley graphs of $\Gamma$ (see Section \ref{percolation}). This beautiful idea is due to Gaboriau \cite{gaboriau-percolation}. Gaboriau and Lyons exploit this idea and give two different proofs of the first step, one using random forests, the other using Bernoulli percolation. They also suggest at the end of their article that the {\em free minimal spanning forest} \cite{lps} could serve as the desired treeable non-hyperfinite subequivalence relation $\mathcal R$. It is this approach that we will present in this paper. Sections \ref{relations} through \ref{subequivalence-bernoulli} are entirely devoted to giving a self-contained proof of this first step. The proof is a combination of ideas and techniques involving probability, ergodic theory, geometric group theory and von Neumann algebras theory.

In the {\bf second step}, one uses Gaboriau's theory of cost \cite{cout} (see also \cite{kechris-miller}). An ergodic treeable non-hyperfinite equivalence relation has cost greater than $1$ by \cite[Th\'eor\`eme ${\rm IV}$.1]{cout}. From the first step, one can then construct an ergodic treeable subequivalence relation $\mathcal R \subset \mathcal R(\Gamma \curvearrowright [0, 1]^\Gamma)$ with cost $\geq 2$. Finally, one applies Hjorth's result \cite{hjorth-cost} in order to get a subequivalence relation of $\mathcal R(\Gamma \curvearrowright [0, 1]^\Gamma)$ induced by a free ergodic pmp action of $\F_2$.

\subsection*{Orbit equivalence theory of nonamenable groups}

As mentioned before, any nonamenable group admits at least two non-OE free ergodic pmp actions \cite{{cw}, {hjorth-T}, {schmidt}}. Over the last few years, the following classes of nonamenable groups have been shown to admit uncountably many non-OE free ergodic pmp actions: property (T) groups (Hjorth \cite{hjorth-T}); nonabelian free groups (Gaboriau and Popa \cite{gaboriau-popa}); weakly rigid groups\footnote{A countable $\Gamma$ is weakly rigid in the sense of Popa if it admits an infinite normal subgroup $\Lambda < \Gamma$ such that the pair $(\Gamma, \Lambda)$ has the relative property (T).} (Popa \cite{popa-cohomology}); nonamenable products of infinite groups (Popa \cite{popasup}, see also  \cite{{monod-shalom}, {ioana04}}); mapping class groups (Kida \cite{kida}). We refer to \cite{{bezuglyi-golodets}, {gefter-golodets}, {zimmer}} for earlier results on this topic.

In his breakthrough paper \cite{ioana07}, Ioana proved that every nonamenable group $\Gamma$ that contains $\F_2$ as a subgroup admits uncountably many non-OE free ergodic pmp actions. As we will see in Section \ref{uncountable}, Ioana's proof goes in two steps that we outline. Regard $\F_2 < \SL_2(\Z)$ as a finite index subgroup and let $\F_2$ act on $\Z^2$ by matrix multiplication. By results of Kazhdan-Margulis \cite{{kazhdan}, {margulis}}, the pair $(\Z^2 \rtimes \F_2, \Z^2)$ has the relative property (T). Write $\alpha : \F_2 \curvearrowright (\mathbf{T}^2, \lambda^2)$ for the corresponding pmp action. The first step (see Theorem \ref{separability}) shows that in every uncountable set of mutually OE actions of $\Gamma$ whose restrictions to $\F_2$ admit $\alpha$ as a quotient, we can find two actions whose restrictions to $\F_2$ are conjugate. The proof is based on a separability argument which uses in a crucial way the fact that the action $\alpha : \F_2 \curvearrowright \mathbf{T}^2$ is {\em rigid} in the sense of Popa \cite{popa2001}. Note that the action $\alpha$ was already successfully used by Gaboriau and Popa \cite{gaboriau-popa} in order to show that the free groups $\F_n$ have a continuum of non-OE actions. The second step consists in using the co-induction technique (see Section \ref{induction}) in order to construct uncountably many actions of $\Gamma$ whose restrictions to $\F_2$ are non-conjugate. Altogether, one obtains uncountably many non-OE actions of $\Gamma$.

Gaboriau and Lyons' result opened up the possibility that the condition ``$\Gamma$ contains $\F_2$" in Ioana's theorem could be replaced by the natural condition ``$\Gamma$ is nonamenable". In order to do so, one had to generalize the second step of Ioana's proof, that is, one needed a more general co-induction technology for group/measurable subgroup rather than group/subgroup. Epstein \cite{epstein} obtained such a construction (see Section $\ref{induction}$). Since the first step of Ioana's proof remains unchanged for $\Gamma$ containing $\F_2$ as a measurable subgroup, Epstein \cite{epstein} obtained the following result.

\begin{theorem}
Every nonamenable group $\Gamma$ admits uncountably many non-OE free ergodic pmp actions.
\end{theorem}

Since then, this result has been generalized in two ways. First, recall that any free ergodic pmp action $\Gamma \curvearrowright (X, \mu)$ gives rise to a finite von Neumann algebra $L^\infty(X) \rtimes \Gamma$ via the  {\em group measure space construction} of Murray and von Neumann (see Section \ref{vonNeumann}). Two free ergodic pmp actions $\Gamma \curvearrowright (X, \mu)$ and $\Lambda \curvearrowright (Y, \nu)$ are W$^*$-{\em equivalent} if the von Neumann algebras $L^\infty(X) \rtimes \Gamma$ and $L^\infty(Y) \rtimes \Lambda$ are $\ast$-isomorphic. Since the group measure space construction only depends on the orbit structure of the action \cite{singer} (see also \cite{feldman-moore2}), it follows that orbit equivalence implies W$^*$-equivalence. Using Popa's concept of {\em rigid inclusion} of von Neumann algebras \cite{popa2001}, Ioana \cite{ioana07} strengthened the previous result by showing that any nonamenable group $\Gamma$ admits a continuum of W$^*$-inequivalent free ergodic pmp actions. Next, given any nonamenable group $\Gamma$, denote by $A_0(\Gamma, X, \mu)$ the standard Borel space of all free mixing pmp actions of $\Gamma$ on $(X, \mu)$ (see \cite{kechris}). On the space $A_0(\Gamma, X, \mu)$, consider the Borel equivalence relation OE defined by $(a, b) \in \operatorname{OE}$ if and only if the actions $a$ and $b$ are orbit equivalent. Epstein, Ioana, Kechris and Tsankov \cite{IKT} proved that OE on the space $A_0(\Gamma, X, \mu)$ cannot be classified by countable structures.

We point out that both Ioana's theorem and Epstein's theorem rely on a separability argument and therefore only provide the {\em existence} of a continuum of non-OE actions for $\Gamma$. What about concrete examples of a continuum of non-OE actions for a given nonamenable group $\Gamma$? Important progress have been made over the recent years. The classes of nonamenable groups for which a concrete uncountable family of non-OE actions is known are the following: non-abelian free groups (Ioana \cite{ioana06}); weakly rigid groups (Popa \cite{popa-cohomology}); nonamenable products of infinite groups (Popa \cite{popasup}); mapping class groups (Kida \cite{kida}). We also refer to Popa and Vaes \cite{popavaes} for further results regarding this question.

\subsection*{Acknowledgments}
I warmly thank Damien Gaboriau for the time we spent discussing the present text and for his numerous valuable suggestions regarding a first draft of this manuscript. I am also grateful to S\'ebastien Martineau for the many thought-provoking conversations regarding percolation theory. I finally thank Vincent Beffara, Adrian Ioana, Russell Lyons, Mikael de la Salle and Stefaan Vaes for their useful comments.

\section{Measure-preserving equivalence relations}\label{relations}

Let $(X, \mu)$ be a nonatomic standard Borel probability space. A {\it countable Borel equivalence relation} $\mathcal{R}$ is an equivalence relation defined on the space $X \times X$ which satisfies:
\begin{enumerate}
\item $\mathcal{R} \subset X \times X$ is a Borel subset.
\item For any $x \in X$, the {\em class} or {\em orbit} of $x$ denoted by $[x]_\mathcal{R} := \{y \in X : (x, y) \in \mathcal{R}\}$ is countable.
\end{enumerate}
We denote by $[\mathcal{R}]$ the {\it full group} of the equivalence relation $\mathcal{R}$, that is, $[\mathcal{R}]$ consists in all Borel isomorphisms $\phi : X \to X$ such that $\graph(\phi) \subset \mathcal{R}$. If $\Gamma$ is a countable group and $(g, x) \to gx$ is a Borel action of $\Gamma$ on $X$, then the {\em orbit equivalence relation} given by
\begin{equation*}
(x, y) \in \mathcal{R}(\Gamma \curvearrowright X) \Longleftrightarrow \exists g \in \Gamma, y = gx
\end{equation*}
is a countable Borel equivalence relation on $X$. By results of Feldman and Moore \cite{feldman-moore}, any countable Borel equivalence relation arises this way. The measure $\mu$ is $\mathcal R$-{\em invariant} if $\phi_*\mu = \mu$, for all $\phi \in [\mathcal R]$. If this is the case, $\mathcal R$ is called a {\em probability measure-preserving} (pmp) equivalence relation on $(X, \mu)$. If $\Gamma \curvearrowright (X, \mu)$ is a pmp action, then $\mathcal{R}(\Gamma \curvearrowright X)$ is a pmp equivalence relation. From now on, we will always assume that $\mathcal R$ is a pmp equivalence relation. Let $\mathcal S$ be a pmp equivalence relation on the nonatomic standard Borel probability space $(Y, \nu)$. We say that $\mathcal R$ and $\mathcal S$ are {\em orbit equivalent} if there exists a pmp Borel isomorphism $\Delta : (X, \mu) \to (Y, \nu)$ such that 
$$(x, y) \in \mathcal R \Longleftrightarrow (\Delta(x), \Delta(y)) \in \mathcal S.$$

For any non-null Borel subset $A \subset X$, define $\mu_A(B) = \mu(B) / \mu(A)$, for all Borel subsets $B \subset A$. Then $(A, \mu_A)$ is a standard Borel probability space. The restricted equivalence relation $\mathcal R \cap (A \times A)$ is simply denoted by $\mathcal R | A$. It is a pmp equivalence relation on $(A, \mu_A)$. The {\em infinite locus} of $\mathcal R$ is the Borel subset 
$$U_\infty := \{x \in X : [x]_\mathcal R \mbox{ is infinite}\}.$$
The restricted equivalence relation $\mathcal R | U_\infty$ is of type ${\rm II_1}$ or aperiodic.\footnote{A pmp equivalence relation $\mathcal R$ is of type ${\rm II_1}$ if almost every $\mathcal R$-class is infinite.} Let $\Gamma \curvearrowright (X, \mu)$ be a {\em free} pmp action of a countable infinite discrete group. Then the orbit equivalence relation $\mathcal R(\Gamma \curvearrowright X)$ induced by the action $\Gamma \curvearrowright X$ is of type ${\rm II_1}$.

For any Borel subset $B \subset X$, define the $\mathcal{R}$-{\it saturation} of $B$ by 
$$[B]_\mathcal{R}  =  \bigcup_{x \in B} [x]_\mathcal{R} =  \{y \in X : \exists x \in B, (x, y) \in \mathcal{R}\}.$$
We have $B \subset [B]_\mathcal{R}$ and $[B]_\mathcal{R}$ is a measurable subset of $X$. We say that $B \subset X$ is $\mathcal{R}$-{\it invariant} if $[B]_\mathcal{R} = B$. The equivalence relation $\mathcal{R}$ is {\it ergodic} if any $\mathcal{R}$-invariant measurable subset $B \subset X$ is null or co-null.  Equivalently, $\mathcal R$ is ergodic if and only if any $[\mathcal R]$-invariant measurable subset $A \subset X$ is null or co-null.

An equivalence relation $\mathcal R$ is {\em hyperfinite} if $\mathcal R = \bigcup_n \mathcal R_n$, where $\mathcal R_n$ is an increasing sequence of {\em finite} subequivalence relations, that is, every orbit of $\mathcal R_n$ is finite. If $\mathcal R$ is hyperfinite, then $\mathcal R | A$ is still hyperfinite for every non-null Borel subset $A \subset X$. Dye \cite{{dye1}, {dye2}} proved there is a unique ergodic hyperfinite ${\rm II_1}$ equivalence relation up to orbit equivalence. It is induced by any ergodic action of $\Z$ on $(X, \mu)$. Ornstein and Weiss \cite{ow} (see also \cite{cfw}) proved that every ergodic pmp action of any infinite amenable group induces the unique ergodic hyperfinite ${\rm II_1}$ equivalence relation.

An ergodic type ${\rm II_1}$ equivalence relation $\mathcal R$ is {\em strongly ergodic} if for every sequence of Borel measurable subsets $A_n \subset X$, we have the following implication: if for all $g \in [\mathcal R]$, we have that $\lim_n \mu(A_n \triangle g A_n) = 0$, then $\lim_n \mu(A_n)(1 - \mu(A_n)) = 0$. A hyperfinite equivalence relation is never strongly ergodic. Let $\Gamma \curvearrowright I$ be any countable infinite group $\Gamma$ acting on a countable set $I$ with infinite orbits and such that for all $g \neq 1_\Gamma$, there are infinitely many $i \in I$ such that $g \cdot i \neq i$. Let $(Y, \nu)$ be any non-trivial probability space and let $(X, \mu) =(Y, \nu)^I$ be the product probability space. The {\em generalized Bernoulli shift} $\Gamma \curvearrowright (Y, \nu)^I$ is defined by $g \cdot (y_i)_{i \in I} = (y_{g^{-1}i})_{i \in I}$. It is a free ergodic pmp action. Moreover, when $\Gamma$ is nonamenable and the action $\Gamma \curvearrowright I$ has amenable stabilizers, the orbit equivalence relation $\mathcal R(\Gamma \curvearrowright Y^I)$ is strongly ergodic. We will use the following characterization of strong ergodicity due to Gaboriau \cite[Proposition 5.2]{gaboriau-icm}.

\begin{prop}\label{strong-ergodicity}
Let $\mathcal R$ be an ergodic type ${\rm II_1}$ equivalence relation on $(X, \mu)$. Then $\mathcal R$ is strongly ergodic if and only if for every increasing sequence $\mathcal R_n$ of subequivalence relations such that $\mathcal R = \bigcup_n \mathcal R_n$, there exist $n \in \N$ and a non-null Borel subset $A \subset X$ such that $\mathcal R_n | A$ is ergodic.
\end{prop}

A pmp {\em graphing} on $(X, \mu)$ is a countable family $\Phi = (\varphi_i)_{i \in I}$ of measure-preserving Borel partial isomorphisms $\varphi_i : A_i \to B_i$. We denote by $\mathcal R_\Phi$ the smallest equivalence relation containing $\{(x, \varphi_i(x)): x \in A_i, i \in I\}$. Then $\mathcal R_\Phi$ is a countable pmp equivalence relation. We say that $\Phi$ {\em generates} the equivalence relation $\mathcal R_\Phi$. The pmp graphing $\Phi$ provides a natural connected graph structure on each class of $\mathcal R$, called the {\em Cayley graph} \cite{cout}. The vertices are the elements of the $\mathcal R$-class and an oriented edge joins two vertices $x$ and $y$ if $x \in A_i$ and $y = \varphi_i(x)$. We denote by $\Phi(x)$ the Cayley graph of $[x]_{\mathcal R}$. A {\em treeing} $\Phi$ is a graphing such that $\mu$-a.s. $\Phi(x)$ is a tree. An equivalence relation $\mathcal R$ is {\em treeable} if there exists a treeing $\Phi$ for which $\mathcal R = \mathcal R_\Phi$. Any hyperfinite equivalence relation is treeable.

The notion of {\em cost} was introduced by Levitt \cite{levitt}. The cost of a pmp graphing $\Phi = (\varphi_i)_{i \in I}$ is defined as $\cost(\Phi, \mu) = \sum_{i \in I} \mu(A_i)$. The cost of a pmp equivalence relation $\mathcal R$ is then defined by
$$\cost(\mathcal R, \mu) = \inf \{\cost(\Phi, \mu) :  \Phi \mbox{ graphing such that } \mathcal R = \mathcal R_\Phi \}.$$
Any ${\rm II_1}$ equivalence relation $\mathcal R$ satisfies $\cost(\mathcal R, \mu) \geq 1$ by \cite{levitt}. Gaboriau proved \cite[Th\'eor\`eme ${\rm IV}.1$]{cout} that when $\mathcal R$ is treeable, $\cost(\mathcal R, \mu) = \cost(\Phi, \mu)$, for every treeing $\Phi$ of $\mathcal R$. In particular when $\mathcal R$ is treeable, $\cost(\mathcal R, \mu) = 1$ if and only if $\mathcal R$ is hyperfinite.

\section{Invariant bond percolation}\label{percolation}

This section is devoted to reviewing a few concepts involving invariant bond percolation on infinite graphs. Further information on this topic may be found in the book \cite{LP} by Lyons and Peres.

\subsection{Graph-theoretic terminology}

Let $\mathcal G = (\mathsf V, \mathsf E)$ be an infinite graph with vertex set $\mathsf V$ and (symmetric) edge set $\mathsf E$. We allow multiple edges and loops. When there is at least one edge joining vertices $u$ and $v$, we say that $u$ and $v$ are adjacent and write $u \sim v$. The {\em degree} $\deg v$ of a vertex $v$ is the number of edges incident with it. A graph is {\em locally finite} if $\deg v < \infty$, for all $v \in \mathsf V$; {\em uniformly bounded} if $\sup_{v \in \mathsf V} \deg v < \infty$; and $d$-regular if $\deg v = d$, for all $v \in \mathsf V$. A connected component of $\mathcal G$ is called a {\em cluster}. A {\em forest} is a graph whose clusters are trees. We will always assume that the graph $\mathcal G$ is locally finite. The automorphism group $\Aut(\mathcal G)$ endowed with the pointwise convergence is locally compact. The graph $\mathcal G$ is {\em transitive} if $\Aut(\mathcal G)$ acts transitively on $\mathsf V$ and {\em unimodular} if $\Aut(\mathcal G)$ is unimodular.

A finite or infinite path $\mathcal P = (\mathbf e_n)_{n \geq 1}$ of edges $\mathbf e_n = [v_n, v_{n + 1}]$ in $\mathcal G$ is {\em self-avoiding} if the map $n \mapsto v_n$ is one-to-one. A {\em simple} cycle is a finite self-avoiding path $\mathcal P = (\mathbf e_1, \dots, \mathbf e_n)$ which is a cycle as well. An {\em infinite} simple cycle is a bi-infinite self-avoiding path $\mathcal P = (\mathbf e_n)_{n \in \mathbf Z}$.

Let $\Gamma$ be a finitely generated group and $S = (s_1, \dots, s_d)$ a finite generating family\footnote{It means that we allow $S$ to contain several copies of the same generator.} for $\Gamma$. Then the (right) {\em Cayley} graph $\mathcal G := \Cay(\Gamma, S)$ is the graph with vertices $\mathsf V := \Gamma$ and edges $\mathsf E := \Gamma \times \{1, \dots, d\}$. The non-oriented edge corresponding to $(v, i)$ will be simply denoted by $[v, vs_i]$. The group $\Gamma$ acts on its Cayley graph by left multiplication. Note that $\Cay(\Gamma, S)$ is an $|S|$-regular transitive unimodular connected graph.

An infinite set of vertices $V$ is {\em end-convergent} if for every finite subset $K \subset \mathcal G$, there is a connected component of $\mathcal G \setminus K$ that contains all but finitely many
vertices of $V$. Two end-convergent sets $V$ and $W$ are equivalent if $V \cup W$ is
end-convergent. An {\em end} of $\mathcal G$ is an equivalence class of end-convergent sets.

\subsection{Bernoulli bond percolation}

In this subsection, we fix an infinite locally finite graph $\mathcal G = (\mathsf V, \mathsf E)$ with $\Gamma < \Aut(\mathcal G)$ a countable discrete subgroup which acts transitively on $\mathsf V$. When $\mathcal G = \Cay(\Gamma, S)$ is the Cayley graph of a finitely generated group $\Gamma$, we regard $\Gamma$ as a discrete subgroup of $\Aut(\mathcal G)$.

We denote by $\{0, 1\}^{\mathsf E}$ the standard Borel space of all subsets of $\mathsf E$, where we identify a subset $A \subset \mathsf E$ with its characteristic function $\mathbf 1_A$. We will regard $\{0, 1\}^{\mathsf E}$ as the Borel space of all subgraphs of $\mathcal G$ with the same set of vertices $\mathsf V$. Observe that $\Gamma$ acts in a Borel way on $\{0, 1\}^{\mathsf E}$ by $(g \cdot \omega)(\mathbf e) = \omega(g^{-1} \mathbf e)$, for all $\mathbf e \in \mathsf E$. Following \cite{{BLPSgafa}, {benjamini-schramm}, {LP}}, a $\Gamma$-invariant {\em bond percolation} $\mathbb P$ on $\mathcal G$ is a $\Gamma$-invariant probability measure on $\{0, 1\}^{\mathsf E}$. The percolation $\mathbb P$ is {\em ergodic} if the pmp action $\Gamma \curvearrowright (\{0, 1\}^\mathsf{E}, \mathbb P)$ is ergodic. We sometimes regard $\omega$ as a $\{0, 1\}^{\mathsf E}$-valued random variable whose law is given by $\mathbb P$. It is customary to denote by $C(\omega; v)$ the cluster of $\omega$ containing the vertex $v$.

For any measurable subset $\mathcal A \subset \{0, 1\}^{\mathsf E}$ and any edge $\mathbf e \in \mathsf E$, denote by $\Pi_{\mathbf e}\mathcal A \subset \{0, 1\}^{\mathsf E}$ the measurable subset $\{\omega \cup \{\mathbf e\} : \omega \in \mathcal A\}$. Likewise denote by $\Pi_{\lnot \mathbf e}\mathcal A \subset \{0, 1\}^{\mathsf E}$ the measurable subset $\{\omega - \{\mathbf e\} : \omega \in \mathcal A\}$. The percolation $\mathbb P$ is {\em insertion tolerant} (resp.\ {\em deletion tolerant}) if for all measurable subset $\mathcal A \subset \{0, 1\}^{\mathsf E}$ such that $\mathbb P[\mathcal A] > 0$ and all  $\mathbf e \in \mathsf E$, we have $\mathbb P[\Pi_{\mathbf e}\mathcal A] > 0$ (resp.\ $\mathbb P[\Pi_{\lnot \mathbf e}\mathcal A] > 0$).

For $p \in [0, 1]$, Bernoulli$(p)$ bond percolation is the product probability measure $\mathbf P_p$ on $\{0, 1\}^{\mathsf E}$ that satisfies $\mathbf P_p [\omega : e \in \omega] = p$. In other words, each edge of $\mathcal G$ is independently kept (or {\em open}) with probability $p$ and removed (or {\em closed}) with probability $1 - p$. The percolation $\mathbf P_p$ is clearly invariant. If the action $\Gamma \curvearrowright \mathsf E$ has infinite orbits, then $\mathbf P_p$ is ergodic. In particular, when $\mathcal G$ is a Cayley graph of an infinite group, $\mathbf P_p$ is ergodic. It is easy to check that $\mathbf P_p$ is both insertion and deletion tolerant for $p \neq 0$ and $1$.

Let $\mathbf P =\Leb^{\mathsf E}$ be the product probability measure on $[0, 1]^{\mathsf E}$ where $\Leb$ denotes the uniform (Lebesgue) measure on $[0, 1]$. An element of $[0, 1]^{\mathsf E}$ gives a colored
graph, with $[0, 1]$ as set of colors. For each
$p \in [0, 1]$, let $\pi_p : [0, 1]^{\mathsf E} \to \{0, 1\}^{\mathsf E}$ be the $\Aut(\mathcal G)$-equivariant map sending $[0, 1]$-colored graphs to
$\{0, 1\}$-colored ones by only keeping the edges colored in $[0, p)$, that is, for every $x \in [0, 1]^{\mathsf E}$,
$$
\pi_p(x)(\mathbf e) = \left\{
\begin{array}{rcl}
1 & \mbox{ if } & x({\mathbf e}) < p \\
0 & \mbox{ if } & x({\mathbf e}) \geq p.
\end{array}
\right.
$$
The {\em standard coupling} is the family $(\pi_p)_{p \in [0, 1]}$. We have that $(\pi_p)_\ast \mathbf P = \mathbf P_p$, for all $p \in [0, 1]$. The event that there exists an infinite cluster in $\pi_p(x)$ is a tail event. Hence, by Kolmogorov's $0, 1$-law, $\mathbf P [\exists \mbox{ an infinite cluster in } \pi_p(x)] = 0 \mbox{ or } 1$. Moreover, for $p \leq q$, the event that $\pi_p(x)$ has an infinite cluster is contained in the event that $\pi_q(x)$ has an infinite cluster. This allows us to define the {\em critical value} $p_c(\mathcal G) \in [0, 1]$ by
$$\mathbf P [\exists \mbox{ an infinite cluster in } \pi_p(x)] = \left\{ 
\begin{array}{rl}
0 & \mbox{ if } p < p_c(\mathcal G) \\
1 & \mbox{ if } p > p_c(\mathcal G). \\
\end{array}
\right.$$
One checks that for all $p \geq p_c(\mathcal G)$, $\mathbf P$-a.s.\ $p_c(\pi_p(x)) = p_c(\mathcal G)/p$.

From now on, assume that the action $\Gamma \curvearrowright \mathsf E$ has infinite orbits, so that the percolation $\mathbf P_p$ is ergodic. Denote by $\mathsf N(\omega)$ the number of infinite clusters of $\omega \in \{0, 1\}^{\mathsf E}$. Since $\mathsf N(\omega)$ is invariant, it follows that $\mathsf N(\omega)$ is a $\mathbf P_p$-a.s.\ constant function, by ergodicity of $\mathbf P_p$. We denote by $\mathsf N_p \in \N \cup \{\infty\}$ its value. Let us prove now that $\mathsf N_p \in \{0, 1, \infty\}$ (see \cite{newman}). Assume that this is not the case, that is, $\mathsf N_p \in \N \setminus \{0, 1\}$. Then there exists a finite path $\mathcal P = (\mathbf e_1, \dots, \mathbf e_n)$ in $\mathcal G$ such that 
$$\mathbf P_p[\mathcal P \mbox{ connects two distinct infinite clusters of } \omega] > 0.$$ 
Denote by $\mathcal{A}$ this last event and let $\mathcal{B} = \Pi_{\mathbf e_1}\circ \cdots \circ \Pi_{\mathbf e_n} (\mathcal{A})$. Since $\mathbf P_p$ is insertion tolerant, $\mathbf P_p[\mathcal{B}] > 0$. Yet, $\mathsf N_p$ takes a strictly smaller value on $\mathcal{B}$ than on $\mathcal A$, which contradicts the fact that $\mathsf N_p$ is a $\mathbf P_p$-a.s.\ constant function.

When $\mathcal G = (\mathsf V, \mathsf E)$ is a connected locally finite unimodular transitive graph, H\"aggstr\"om and Peres \cite{haggstrom-peres} showed there is {\em monotonicity of uniqueness}: for all $0 \leq p_1 < p_2 \leq 1$, 
$$
\begin{array}{rl}
\mbox{if} & \mathbf P[\exists \mbox{ a unique infinite cluster in } \pi_{p_1}(x)] = 1 \\
\mbox{then} & \mathbf P[\exists \mbox{ a unique infinite cluster in } \pi_{p_2}(x)] = 1. 
\end{array}
$$
This explains why the {\em uniqueness phase} is an interval and allows us to define
$$
p_u(\mathcal G) =  \inf \{p \in [0, 1] : \mbox{ there is a unique infinite cluster for } \mathbf P_p\}.
$$
We have $p_c(\mathcal G) \leq p_u(\mathcal G)$. Stronger still, H\"aggstr\"om and Peres \cite{haggstrom-peres} proved that after $p_c(\mathcal G)$, there is no spontaneous generation of infinite clusters, ``all infinite clusters are born simultaneously":

\begin{theo}\label{hp}
Let $\mathcal G = (\mathsf V, \mathsf E)$ be a connected locally finite unimodular transitive graph. The number $\mathsf N_p$ of infinite clusters in $\pi_p(x)$ is a $\mathbf P$-a.s.\ constant function and we have
$$
\mathsf N_p = \left\{ 
\begin{array}{rcl}
0 & \mbox{for} & p \in [0, p_c(\mathcal G)) \\
\infty & \mbox{for} & p \in (p_c(\mathcal G), p_u(\mathcal G)) \\
1 & \mbox{for} & p \in (p_u(\mathcal G), 1].
\end{array}
\right.
$$
\begin{itemize}
\item Moreover, for all $p_1 < p_2$, when $\mathbf P$-a.s.\ $\pi_{p_1}(x)$ produces at least one infinite cluster, $\mathbf P$-a.s.\ every infinite cluster of $\pi_{p_2}(x)$ contains at least one infinite cluster of $\pi_{p_1}(x)$.
\item If $\mathbf P$-a.s.\ $\pi_p(x)$ produces infinitely many infinite clusters, then $\mathbf P$-a.s.\ all infinite clusters of $\pi_p(x)$ have uncountably many ends.
\item When $p < 1$, if $\mathbf P$-a.s.\ $\pi_p(x)$ produces only one infinite cluster, then $\mathbf P$-a.s.\ the unique infinite cluster of $\pi_p(x)$ has only one end.
\end{itemize}
\end{theo}

Lyons and Schramm \cite{lyons-schramm} showed that when Bernoulli$(p)$ bond percolation produces a.s.\ at least one infinite cluster, then its infinite clusters are {\em indistinguishable} in the following sense. Consider the Borel subset 
$$\mathfrak C_\infty = \left\{(\omega, C) \in 2^{\mathsf E} \times 2^{\mathsf V} : C \mbox{ is an infinite cluster of } \omega \right\}.$$
Observe that $\mathfrak C_\infty$ is invariant under the diagonal action of $\Gamma$. A $\Gamma$-invariant bond percolation $\mathbb P$ on $\mathcal G$ has {\em indistinguishable infinite clusters} if for every $\Gamma$-invariant Borel subset $\mathcal A \subset \mathfrak C_\infty$, $\mathbb P$-a.s.\ either for all infinite clusters $C$ of $\omega$, we have $(\omega, C) \in \mathcal A$, or for all infinite clusters $C$ of $\omega$, we have $(\omega, C) \in \mathfrak C_\infty \setminus \mathcal A$. Observe that when $\mathbb P$ is moreover ergodic, we can permute ``$\mathbb P$-a.s." with ``or". The following result is \cite[Theorem 3.3]{lyons-schramm}.

\begin{theo}[Clusters indistinguishability]\label{indistinguishability}
Let $\mathcal G = (\mathsf V, \mathsf E)$ be a unimodular transitive graph. Any $\Gamma$-invariant insertion-tolerant bond percolation on $\mathcal G$ has indistinguishable infinite clusters.
\end{theo}

\subsection{From percolation to equivalence relations}
Let $\Gamma$ be a finitely generated infinite group and $S = (s_1, \dots, s_d)$ a finite generating family for $\Gamma$. Set $\mathcal G = \Cay(\Gamma, S)$ that we also denote $\mathcal G = (\mathsf V, \mathsf E)$. Let $\Gamma \curvearrowright (X, \mu)$ be a free ergodic pmp action and denote by $\mathcal S := \mathcal R(\Gamma \curvearrowright X)$ the induced orbit equivalence relation. Let $\pi : X \to \{0, 1\}^{\mathsf E}$ be a $\Gamma$-equivariant Borel map. Then the push-forward measure $\pi_\ast\mu$ is a $\Gamma$-invariant bond percolation on $\mathcal G$. The following definition is due to Gaboriau \cite{gaboriau-percolation}. 
\begin{defi}
The {\em cluster} subequivalence relation $\mathcal R^{\cl}_\pi \subset \mathcal S$ is defined by
$$
(x, y) \in \mathcal R^{\cl}_\pi \Longleftrightarrow \left\{
\begin{array}{l}
\mbox{there exists } g \in \Gamma, y = g^{-1} x \\
1_\Gamma \mbox{ and } g \mbox{ are in the same cluster of } \pi(x). 
\end{array}
\right.
$$ 
\end{defi}
Denote by $\mathbf e_i$ the edge $[1_\Gamma, s_i]$. Define the Borel set $X_i = \{x \in X : \pi(x)(\mathbf e_i) = 1\}$ and partial Borel isomorphisms $\varphi_i = s_i^{-1} : X_i \to s_i^{-1}(X_i)$. Then the family $\Phi = (\varphi_1, \dots, \varphi_d)$ is a pmp graphing which generates $\mathcal R^{\cl}_\pi$ and $\Phi(x) \simeq C(\pi(x); 1_\Gamma)$, for $\mu$-almost every $x \in X$. Denote by $U_\infty^\pi$ the infinite locus of $\mathcal R_\pi^{\cl}$, that is,
$$U_\infty^\pi = \{x \in X : C(\pi(x), 1_\Gamma) \mbox{ is infinite}\}.$$
Assume now that $\mu$-a.s.\ $\pi(x)$ produces at least one infinite cluster. Then $\mu(U_\infty^\pi) > 0$ and $\mathcal R_\pi^{\cl} | U_\infty^\pi$ is a type ${\rm II_1}$ equivalence relation. Moreover, on $U_\infty^\pi$, each $\mathcal S$-class splits into $\mathcal R_\pi^{\cl}$-classes which are in one-to-one correspondence with the infinite clusters of $\pi(x)$. It follows in particular that when $\mu$-a.s.\ $\pi(x)$ produces exactly one infinite cluster, the orbit and the cluster equivalence relations do coincide on the infinite locus, that is, $\mathcal R_\pi^{\cl} | U_\infty^\pi = \mathcal S | U_\infty^\pi$. The following  observation is due to Gaboriau and Lyons \cite{gaboriau-lyons}.

\begin{prop}[Indistinguishability vs. ergodicity]
The percolation $\pi_\ast\mu$ has indistinguishable infinite clusters if and only if the equivalence relation $\mathcal R_\pi^{\cl} | U_\infty^\pi$ is ergodic.
\end{prop}

Consider now Bernoulli$(p)$ bond percolation through the standard coupling $(\pi_p)_{p \in [0, 1]}$. Observe that since the action $\Gamma \curvearrowright \mathsf E$ is free, the free pmp action $\Gamma \curvearrowright ([0, 1]^{\mathsf E}, \mathbf P)$ is conjugate to the plain Bernoulli shift $\Gamma \curvearrowright ([0, 1], \Leb)^{\Gamma}$. Let $\mathcal S$ be the corresponding orbit equivalence relation. Simply denote by $\mathcal R_p$ the cluster equivalence relation $\mathcal R^{\cl}_{\pi_p}$. The family $(\mathcal R_p)_{p \in [0, 1]}$ is increasing. Moreover $\mathcal R_q = \bigcup_{p < q} \mathcal R_p$ and $\mathcal R_1 = \mathcal S$.
\begin{itemize}
\item For $p < p_c(\mathcal G)$, $\mathbf P$-almost every orbit of $\mathcal R_p$ is finite, that is, $\mathcal R_p$ is a type $\mathrm I$ equivalence relation. It follows in particular that $\mathcal R_{p_c(\mathcal G)}$ is {\em hyperfinite}.

\item For $p > p_c(\mathcal G)$, denote by $U_\infty^p$ the (non-null) infinite locus of $\mathcal R_p$. If $\mathbf P$-a.s.\ $\pi_p(x)$ produces infinitely many infinite clusters, $\mathcal R_p | U_\infty^p$ has infinite index in $\mathcal S | U_\infty^p$.
\end{itemize}

It is straightforward to see that {\em clusters indistinguishability} implies {\em simultaneous uniqueness}. Indeed, simultaneous uniqueness amounts to saying that for all $p_1 < p_2$ such that $\mathbf P[U_\infty^{p_1}] > 0$, the $\mathcal R_{p_2} | U_\infty^{p_2}$-saturation of $U_\infty^{p_1}$ is equal to $U_\infty^{p_2}$. This is clear since $\mathcal R_{p_2} | U_\infty^{p_2}$ is ergodic by clusters indistinguishability.

\section{The non-uniqueness phase in Bernoulli percolation}

A famous conjecture by Benjamini and Schramm \cite[Conjecture 6]{benjamini-schramm} is that if a transitive graph $\mathcal G$ with finite degree is nonamenable, then $p_c(\mathcal G) < p_u(\mathcal G)$. This section is devoted to presenting a partial answer to this question, due to Pak and Smirnova-Nagnibeda \cite{smirnova}: for any nonamenable finitely-generated group $\Gamma$, there exists a finite  generating family $S$ such that the Cayley graph $\mathcal G := \Cay(\Gamma, S)$ has a non-uniqueness phase, that is, for which $p_c(\mathcal G) < p_u(\mathcal G)$.

Let $\mathcal G = \Cay(\Gamma, S)$ be a Cayley graph of an infinite finitely generated group $\Gamma$ with respect to a finite generating family $S = (s_1, \dots, s_d)$. Recall that the vertex set $\mathsf V$ is $\Gamma$ and the edge set $\mathsf E$ is $\{[g, gs_i] : g \in \Gamma, 1 \leq i \leq d\}$. For a non-empty finite subset $F \subset \mathsf V$, let $\partial_{\mathsf E} F$ be the set of edges which have exactly one endpoint in $F$. Define the {\em edge-isoperimetric} constant of $\mathcal G$ by
$$\iota_{\mathsf E}(\mathcal G) := \inf \left \{ \frac{|\partial_{\mathsf E} F|}{|F|} : \emptyset \neq F \subset \mathsf V \mbox{ finite subset}\right \}.$$
A graph $\mathcal G$ is {\em edge-amenable} if $\iota_{\mathsf E}(\mathcal G) = 0$. A finitely generated group $\Gamma$ is amenable if for some (or equivalently for every) finite generating family $S$, the  Cayley graph $\Cay(\Gamma, S)$ is edge-amenable. The first result of this section is due to Benjamini and Schramm \cite[Theorem 2]{benjamini-schramm}.

\begin{theo}[Upper bound for $p_c$]\label{upper}
Let $\mathcal{G} = \Cay(\Gamma, S)$. Then
$$p_c(\mathcal{G}) \leq \frac{1}{\iota_{\mathsf E}(\mathcal{G}) + 1}.$$
\end{theo}
\begin{proof}
Fix $p > \frac{1}{\iota_{\mathsf E}(\mathcal{G}) + 1}$ and let $\mathbf P_p$ be the corresponding Bernoulli$(p)$ percolation on $\mathcal G$. Fix $v \in \mathsf V$. Let $(\mathbf e_i)_{i \geq 1}$ be an ordering of $\mathsf E$ so that $\mathbf e_1$ is incident with $v$. Let $\omega \in \{0, 1\}^{\mathsf E}$ be a configuration. We explore the open cluster $C(\omega; v)$ by looking at the following inductive procedure.

Let $E_1 = \{\mathbf e_1\}$, $V_1 = \{ v \}$ and $X_1(\omega) = \omega(\mathbf e_1)$. Assume $E_k$ and $V_k$ have been defined. Denote by $V_{k + 1}$ the set $\{v\} \cup \{ \mbox{endpoints of open edges in } E_k\}$. Let $n_{k + 1}$ be the least integer $n$ such that the edge $\mathbf e_n \in \mathsf E \setminus E_k$ has exactly one endpoint in $V_{k + 1}$, if any. 
\begin{enumerate}
\item [(a)] If there are none, then stop. Denote by $k := n(\omega)$ the stopping time. In that case, the open cluster $C(\omega; v)$ containing $v$ is finite. Then set $\ell_k = \sup \{n_j : 1 \leq j \leq k\}$ and $X_{k + i}(\omega) = \omega(\mathbf e_{\ell_k + i})$, for all $i \geq 1$.
\item [(b)] Otherwise, let $E_{k + 1} = E_k \cup \{\mathbf e_{n_{k + 1}}\}$ and $X_{k + 1}(\omega) = \omega(\mathbf e_{n_{k + 1}})$. 
\end{enumerate}
If the procedure never ends, then the open cluster $C(\omega; v)$ is infinite.

\begin{claim}
$(X_n)_{n \geq 1}$ is an infinite sequence of i.i.d.\ $\{0, 1\}$-valued Bernoulli$(p)$ random variables.
\end{claim}

It suffices to show that for all $k \geq 1$ and all $\varepsilon_1, \dots, \varepsilon_k \in \{0, 1\}$, we have
\begin{equation}\label{independence}
\mathbf P_p [X_{k + 1} = 1 | X_1 = \varepsilon_1, \dots, X_k = \varepsilon_k] = p.
\end{equation}
Denote by $\mathcal A = \{\omega : X_1(\omega) = \varepsilon_1, \dots, X_k(\omega) = \varepsilon_k\}$, $\mathcal A_i = \mathcal A \cap \{\omega : n(\omega) = i\}$, for $1 \leq i \leq k$, and $\mathcal A_{k + 1} = \mathcal A \cap \{\omega : n(\omega) \geq k + 1\}$. For $i \leq k$, there are $k +  1$ fixed distinct edges $\mathbf f_1 = \mathbf e_{n_1}, \dots, \mathbf f_i = \mathbf e_{n_i}, \mathbf f_{i + 1} = \mathbf e_{\ell_i + 1}, \dots, \mathbf f_{k + 1} = \mathbf e_{\ell_i + k + 1 - i}$, with $\ell_i = \sup \{n_j : 1 \leq j \leq i\}$, such that $\mathcal A_i = \{\omega : \omega(\mathbf f_1) = \varepsilon_1, \dots, \omega(\mathbf f_k) = \varepsilon_k\}$. We moreover have $\mathbf P_p [X_{k + 1} = 1 | \mathcal A_i] = \mathbf P_p [\omega(\mathbf f_{k + 1}) = 1 | \mathcal A_i]$. Since the edges $\mathbf f_1, \dots, \mathbf f_{k + 1}$ are distinct, the random variables $\omega(\mathbf f_1), \dots, \omega(\mathbf f_{k + 1})$ are independent. It follows that $\mathbf P_p [\omega(\mathbf f_{k + 1}) = 1 | \mathcal A_i] = p$. Likewise, for $i = k + 1$, there are $k + 1$ fixed distinct edges $\mathbf e_{n_1}, \dots, \mathbf e_{n_{k + 1}}$ such that $\mathcal A_{k + 1} = \{\omega : \omega(\mathbf e_{n_1}) = \varepsilon_1, \dots, \omega(\mathbf e_{n_k}) = \varepsilon_k \}$. We moreover have $\mathbf P_p [X_{k + 1} = 1 | \mathcal A_{k + 1}] = \mathbf P_p [\omega(\mathbf e_{n_{k + 1}}) = 1 | \mathcal A_{k + 1}]$. Since the edges $\mathbf e_{n_1}, \dots, \mathbf e_{n_{k + 1}}$ are distinct, the random variables $\omega(\mathbf e_{n_1}), \dots, \omega(\mathbf e_{n_{k + 1}})$ are independent. It follows that $\mathbf P_p [\omega(\mathbf e_{n_{k + 1}}) = 1 | \mathcal A_{k + 1}] = p$. Since the event $\mathcal A$ is equal to the disjoint union of the events $\mathcal A_1, \dots, \mathcal A_{k + 1}$, we get Equation $(\ref{independence})$, which finishes the proof of the claim.

By the strong law of large numbers, we get
$$\mathbf P_p\left[\frac1n \sum_{k = 1}^n X_k(\omega) > \frac{1}{\iota_{\mathsf E}(\mathcal{G}) + 1}, \forall n \geq 1\right] > 0.$$
We denote by $\mathcal A$ this last event. We show that $C(\omega; v)$ must be infinite on the event $\mathcal{A}$. Assume that $C(\omega; v)$ is finite. Simply denote $n = n(\omega)$ and let $E_n$ be the last set of selected edges according to (a). Let $m = |C(\omega; v)|$. We have that $E_n$ contains $\partial_\mathsf{E} C(\omega; v)$ (for which all edges are closed) and a spanning tree of $C(\omega; v)$ with $m - 1$ open edges. Thus we have $n \geq |\partial_{\mathsf E} C(\omega; v)| + m - 1$ and $\sum_{k = 1}^n X_k(\omega) = m - 1$, so that  
$$\frac1n \sum_{k = 1}^n X_k(\omega) = \frac{m - 1}{n} \leq \frac{m - 1}{ |\partial_{\mathsf E} C(\omega; v)| + m - 1} \leq \frac{1}{\frac{|\partial_{\mathsf E} C(\omega; v)|}{|C(\omega; v)|} + 1} \leq \frac{1}{\iota_{\mathsf E}(\mathcal G) + 1}.$$
It follows that $C(\omega; v)$ is infinite on the event $\mathcal A$ and thus 
$$\mathbf P_p [C(\omega; v) \mbox{ is infinite}] > 0.$$
Therefore $p > p_c (\mathcal G)$, which finishes the proof.
\end{proof}

Let $\mathcal G = \Cay(\Gamma, S)$, where $S = (s_1, \dots, s_d)$. Let $P : \ell^2(\Gamma) \to \ell^2(\Gamma)$ be the corresponding {\em simple random walk} operator: for all $f \in \ell^2(\Gamma)$, 
$$(Pf)(g) = \frac1d \sum_{i = 1}^d f(g s_i).$$
It is easy to see that as a bounded operator on $\ell^2(\Gamma)$, we have $P = P^*$ and $\|P\|_\infty \leq 1$ (where $\| \cdot \|_\infty$ is the operator norm). Fix an orientation of the edges. Define the {\em differential} operator $\partial : \ell^2(\Gamma) \to \ell^2(\mathsf E)$ by $(\partial f)(\mathbf e) = f(\mathbf e_+) - f(\mathbf e_-)$. The combinatorial Laplacian is then defined as the positive self-adjoint operator $\Delta  = \partial^*\partial$. A straightforward computation gives $\Delta = d (1 - P)$. The {\em spectral radius} of the graph $\mathcal G$ is defined as $\rho(\mathcal G) := \|P\|_\infty$.

\begin{prop}[\cite{mohar}]\label{inequality}
Let $\mathcal G = \Cay(\Gamma, S)$, where $S = (s_1, \dots, s_d)$. Then
$$\iota_{\mathsf E}(\mathcal G) \geq d(1 - \rho(\mathcal G)).$$
\end{prop}

\begin{proof}
Let $F \subset \mathsf V$ be a nonempty finite subset. Let $f = \mathbf 1_F$. We have 
$$|\partial_{\mathsf E} F| = \langle \Delta f, f\rangle = d \langle (1 - P) f, f\rangle \geq d (1 - \rho(\mathcal G)) \|f\|^2 = d (1 - \rho(\mathcal G)) |F|,$$
and the Proposition follows.
\end{proof}

Choose a vertex $v \in \mathsf V$ (e.g.\ $v = 1_\Gamma$) and denote by $a_n(\mathcal G)$ the number of simple cycles of length $n$ in $\mathcal G$ that contain $v$. Let
$$\gamma(\mathcal{G}) := \limsup_n a_n(\mathcal{G})^{1/n}.$$
Denote by $(\langle X_n\rangle, \mathbf P_v)$ the simple random walk on $\mathcal G$ starting at $v$. Recall that $\rho(\mathcal G) = \limsup_n (\mathbf P_v[X_n = v])^{1/n}$. Any simple cycle of length $n$ that contains $v$ defines a way for the simple random walk starting at $v$ to return to $v$ at time $n$. That event has probability $1/d^n$. Therefore $\mathbf P_v[X_n = v] \geq a_n(\mathcal G)/d^n$, which shows that $\gamma(\mathcal G) \leq d \rho(\mathcal G)$. The next theorem, due to Schramm, is an improvement of an earlier result of Benjamini and Schramm \cite[Theorem 4]{benjamini-schramm}. The proof we give here is borrowed from Lyons \cite[Theorem 3.9]{lyons2000}.
\begin{theo}[Lower bound for $p_u$]\label{lower}
Let $\mathcal{G} = \Cay(\Gamma, S)$. Then
$$\frac{1}{\gamma(\mathcal G)} \leq p_u(\mathcal G).$$
\end{theo}
\begin{proof}
Let $1 > p > p_u(\mathcal{G}) \geq p_c(\mathcal{G})$. Since $p > p_u(\mathcal{G})$, we know that $\mathbf P_p$-a.s.\ the open subgraph $\omega$ contains a unique infinite cluster $C(\omega)$ which has only one end. We start by proving the following.

\begin{claim}[\cite{lps}]
Let $G$ be a graph of bounded degree that does not contain an infinite simple cycle. Then $p_c(G) = 1$.
\end{claim}

By repeated applications of Menger's Theorem\footnote{For any vertex $v$ in an infinite graph $\mathcal G$, the maximum number of paths from $v$ to $\infty$ that are pairwise disjoint (except at $v$) is equal to the minimum cardinality of a set $W$ of
vertices such that $W$ is disjoint from $v$, but every path from $v$ to $\infty$ passes through $W$.} we see that if $v$ is a vertex in $G$, then there are infinitely many vertices $v_n$ such that $v$ is in a finite cluster of $G \setminus \{v_n\}$. Since $G$ has bounded degree, it follows that $p_c(G) = 1$, which finishes the proof of the claim.

We get that $\mathbf P_p$-a.s.\ $\omega$ contains an infinite simple cycle. Otherwise, the claim would imply that with $\mathbf P_p$-positive probability, $p_c(\omega) = 1$. This contradicts the fact that $\mathbf P_p$-a.s.\ $p_c(\omega) = p_c(\mathcal G)/ p < 1$.

Denote by $\mathcal A \subset \{0, 1\}^{\mathsf E}$ the event that there is an infinite simple cycle in the $p$-open cluster $C(\omega)$ containing $v$. We may regard such an infinite simple cycle as the union of two disjoint infinite simple rays starting at $v$. We have proven that $\mathbf P_p[\mathcal A] > 0$. Since $C(\omega)$ has only one end, these two paths may be connected by paths in $\omega$ that
stay outside arbitrarily large balls. In particular, there are an infinite number of simple cycles in $\omega \in \mathcal A$ through the vertex $v$. The expected number of such simple cycles must be infinite, whence we obtain in particular $\sum_n a_n(\mathcal{G}) p^n = \infty$. Thus $p  > \gamma(\mathcal G)^{-1}$, which finishes the proof.
\end{proof}

\begin{coro}\label{nonuniqueness}
Let $\mathcal G = \Cay(\Gamma, S)$. Assume that $\rho(\mathcal G) \leq 1/2$. Then $p_c(\mathcal G) < p_u(\mathcal G)$.
\end{coro}

\begin{proof}
Using Proposition \ref{inequality}, Theorems \ref{upper} and \ref{lower}, we have 
$$p_c(\mathcal G) \leq  \frac{1}{\iota_{\mathsf E}(\mathcal{G}) + 1} <  \frac{1}{\iota_{\mathsf E}(\mathcal{G})} \leq  \frac{1}{d(1 - \rho(\mathcal G))} \leq  \frac{1}{d \rho(\mathcal G)} \leq \frac{1}{\gamma(\mathcal G)} \leq p_u(\mathcal G).$$
\end{proof}

We finally state and prove the result of Pak and Smirnova-Nagnibeda \cite{smirnova}.

\begin{coro}\label{pak-smirnova}
Let $\Gamma$ be a finitely generated nonamenable group. Then there exists a generating family $S$ of $\Gamma$ such that $p_c(\Cay(\Gamma, S)) < p_u(\Cay(\Gamma, S))$.
\end{coro}

\begin{proof}
Let $S$ be a finite generating family for $\Gamma$ such that $1_\Gamma \in S$ and let $\mathcal G = \Cay(\Gamma, S)$. For $k \geq 1$, define the $k$-fold family $S^{[k]}$.
The group $\Gamma$ may be regarded as generated by $S^{[k]}$. Let $\mathcal G^{[k]} = \Cay(\Gamma, S^{[k]})$. If $P$ denotes the random walk operator on the graph $\mathcal G$, then $P^k$ is the random walk operator of $\mathcal G^{[k]}$. Thus 
$$\rho(\mathcal G^{[k]}) = \|P^k\|_\infty \leq \|P\|_\infty^k = \rho(\mathcal G)^k.$$
Since $\Gamma$ is nonamenable, $\rho(\mathcal G) < 1$ by Kesten's result \cite{kesten}. Let $k$ be a large enough integer so that $\rho(\mathcal G)^k \leq 1/2$. We finally get $\rho(\mathcal G^{[k]}) \leq 1/2$. By Corollary \ref{nonuniqueness}, the finite generating family $S^{[k]}$ does the job.
\end{proof}

\section{Minimal spanning forests and applications}

\subsection{Minimal spanning forests}

We first review results due to Lyons, Peres and Schramm \cite{lps} regarding {\em minimal spanning forests} on infinite connected graphs and their relation to Bernoulli percolation.

Let $\mathcal G = \Cay(\Gamma, S)$ be a Cayley graph of an infinite finitely generated group $\Gamma$ with respect to a finite generating family $S$. As usual, denote by $\mathsf V$ the vertex set and by $\mathsf E$ the edge set. Denote by $\forest(\mathcal G) \subset \{0, 1\}^{\mathsf E}$ the Borel subset of all forests of $\mathcal G$. A {\em random forest} is an invariant bond percolation supported on $\forest(\mathcal G)$. We endow the Borel space $[0, 1]^{\mathsf E}$ with the product probability measure $\mathbf P = \Leb^{\mathsf E}$. Given $x \in [0, 1]^{\mathsf E}$ an injective labeling of the edges, let $\mathsf{FMSF}(x)$ be the set of edges $\mathbf e \in \mathsf E$ such that in every simple cycle in $\mathcal G$ containing $\mathbf e$, there exists at least one edge $\mathbf e' \neq \mathbf e$ with $x(\mathbf e') > x(\mathbf e)$. The $\Aut(\mathcal G)$-equivariant map $\mathsf{FMSF} : [0, 1]^{\mathsf E} \to \{0, 1\}^{\mathsf E}$ (or simply its law) is called the {\em free minimal spanning forest} on $\mathcal G$. Observe that if $\mathcal G$ is a tree, then $\mathbf P$-a.s.\ $\mathsf{FMSF}(x) = \mathcal G$.

An {\em extended simple cycle} in $\mathcal G$ is either a simple cycle in $\mathcal G$ or an infinite simple cycle in $\mathcal G$. Given $x \in [0, 1]^{\mathsf E}$ an injective labeling of the edges, let $\mathsf{WMSF}(x)$ be the set of edges $\mathbf e \in \mathsf E$ such that in every extended simple cycle in $\mathcal G$ containing $\mathbf e$, there exists at least one edge $\mathbf e' \neq \mathbf e$ with $x(\mathbf e') > x(\mathbf e)$. Equivalently, $\mathsf{WMSF}(x)$ consists of those edges $\mathbf e$ such that there is a finite set $W \subset \mathsf V$ where $\mathbf e$ is the least edge joining $W$ to $\mathsf V \setminus W$. The $\Aut(\mathcal G)$-equivariant map $\mathsf{WMSF} : [0, 1]^{\mathsf E} \to \{0, 1\}^{\mathsf E}$ (or simply its law) is called the {\em wired minimal spanning forest} on $\mathcal G$. Observe that if $\mathcal G$ is a tree with one end, then $\mathbf P$-a.s.\ $\mathsf{WMSF}(x) = \mathcal G$.

It is clear that $\mathsf{WMSF}(x) \subset \mathsf{FMSF}(x)$. Moreover, $\mathsf{WMSF}(x)$ and $\mathsf{FMSF}(x)$ are indeed forests since in every simple cycle in $\mathcal G$, the edge $\mathbf e$ with maximum label $x(\mathbf e)$ is contained neither in $\mathsf{WMSF}(x)$ nor in $\mathsf{FMSF}(x)$. Moreover, all the clusters of $\mathsf{WMSF}(x)$ and $\mathsf{FMSF}(x)$ are infinite since the least edge joining every finite vertex set to its complement belongs to both forests.

Define
$$f(x, \mathbf e) := \inf_{\mathcal P} \max \{x(\mathbf e') : \mathbf e' \in \mathcal P, \mathbf e' \neq \mathbf e\},$$
where the infimum is over simple cycles $\mathcal P$ that contain the edge $\mathbf e$. If there are none, the infimum is defined to be $\infty$. It follows that $\mathsf{FMSF}(x) = \{\mathbf e \in \mathsf E : x(\mathbf e) \leq f(x, \mathbf e)\}$. Likewise, define
$$w(x, \mathbf e) := \inf_{\mathcal P} \sup \{x(\mathbf e') : \mathbf e' \in \mathcal P, \mathbf e' \neq \mathbf e\},$$
where the infimum is over extended simple cycles $\mathcal P$ in $\mathcal G$ that contain the edge $\mathbf e$. If there are none, the infimum is defined to be $\infty$. It follows that $\{\mathbf e \in \mathsf E : x(\mathbf e) < w(x, \mathbf e)\} \subset \mathsf{WMSF}(x) \subset \{\mathbf e \in \mathsf E : x(\mathbf e) \leq w(x, \mathbf e)\}$. Since $x(\mathbf e)$ and $w(x, \mathbf e)$ are independent random variables and $x(\mathbf e)$ is uniformly distributed, we get $\mathbf P$-a.s.\ 
$$\mathsf{WMSF}(x) = \{\mathbf e \in \mathsf E : x(\mathbf e) < w(x,\mathbf e)\} = \{\mathbf e \in \mathsf E : x(\mathbf e) \leq w(x,\mathbf e)\}.$$ It is clear that $w(x, \mathbf e) \leq f(x, \mathbf e)$, for all $\mathbf e \in \mathsf E$. The following is \cite[Proposition 6]{lps}.

\begin{prop}\label{forests}
Let $\mathcal G = \Cay(\Gamma, S)$. Then $\mathsf{WMSF} \neq \mathsf{FMSF}$ if and only if $p_c(\mathcal G) < p_u(\mathcal G)$.
\end{prop}

\begin{proof}
We will use the standard coupling $\pi_p : ([0, 1]^{\mathsf E}, \mathbf P) \to (\{0, 1\}^{\mathsf E}, \mathbf P_p)$ as defined previously. Since $\mathsf{WMSF}(x) \subset \mathsf{FMSF}(x)$ and $\mathsf E$ is countable, it follows that $\mathsf{WMSF} \neq \mathsf{FMSF}$ if and only if there exists $\mathbf e \in \mathsf E$ such that $ \mathbf P[w(x, \mathbf e) < x(\mathbf e) \leq f(x, \mathbf e)] > 0$. Recall that $x(\mathbf e)$ is independent from the random variables $w(x, \mathbf e)$ and $f(x,\mathbf e)$, and $x(\mathbf e)$ is uniformly distributed. Therefore $\mathsf{WMSF} \neq \mathsf{FMSF}$ if and only if there exist $\mathbf e \in \mathsf E$  and $p_1 < p_2$ such that $ \mathbf P[w(x, \mathbf e) \leq p_1 < p_2 \leq f(x, \mathbf e)] > 0$.

Assume that $p_c(\mathcal G) < p_u(\mathcal G)$. Let $p_c(\mathcal G) < p_1 < p_2 < p_u(\mathcal G)$. Using Theorem \ref{hp}, we know that $\mathbf P$-a.s.\ $\pi_{p_2}(x)$ has at least two distinct infinite clusters and each of these clusters contains an infinite cluster of $\pi_{p_1}(x)$. Therefore there exists a simple path $\mathcal P = (\mathbf e_1, \dots, \mathbf e_n)$ of minimal length $n$ in $\mathcal G$, where $\mathbf e_i = [v_i, v_{i + 1}]$, such that with $\mathbf P$-positive probability the following hold:
\begin{enumerate}
\item $\mathcal P$ connects two distinct infinite clusters of $\pi_{p_1}(x)$.
\item The clusters $C(\pi_{p_2}(x); v_1)$ and $C(\pi_{p_2}(x); v_{n + 1})$ are infinite and distinct. 
\end{enumerate}
Using the standard coupling and since $\mathbf P_{p_1}$ and $\mathbf P_{p_2}$ are both insertion and deletion tolerant, the minimal length of $\mathcal P$ has to be $1$. In other words, there exists an edge $\mathbf e \in \mathsf E$ such that with $\mathbf P$-positive probability, the two endpoints of $\mathbf e$ are in distinct infinite clusters of $\pi_{p_i}(x)$, for $i = 1, 2$. We get $\mathbf P[w(x, \mathbf e) \leq p_1 < p_2 \leq f(x, \mathbf e)] > 0$, whence $\mathsf{WMSF} \neq \mathsf{FMSF}$.

Conversely, assume that $\mathsf{WMSF} \neq \mathsf{FMSF}$. In particular, there exist $\mathbf e \in \mathsf E$  and $p$ such that $\mathbf P[w(x, \mathbf e) < p \leq f(x, \mathbf e)] > 0$. 
Then $\mathbf P[w(x, \mathbf e) < p \leq f(x, \mathbf e) \mbox{ and } p \leq x(\mathbf e)] > 0$. It follows that with $\mathbf P$-positive probability, $\pi_p(x)$ has at least two distinct infinite clusters, whence $p_c(\mathcal G) < p_u(\mathcal G)$.
\end{proof}

\subsection{Cluster equivalence relations of $\mathsf{MSF}$}

We denote by $\mathcal R_{\mathsf{WMSF}}$ and $\mathcal R_{\mathsf{FMSF}}$ the cluster equivalence relations associated to both minimal spanning forests on $\mathcal G = \Cay(\Gamma, S)$. Both of them are of type ${\rm II_1}$ and the treeing of $\mathcal R_{\mathsf{WMSF}}$ is a subtreeing of $\mathcal R_{\mathsf{FMSF}}$, that is, $\mathcal R_{\mathsf{WMSF}} \subset \mathcal R_{\mathsf{FMSF}}$. Lyons, Peres and Schramm proved that $\mathbf P$-a.s.\ every tree of $\mathsf{WMSF}(x)$ has exactly one end (see \cite[Theorem 3.12]{lps}). In other words, $\mathcal R_{\mathsf{WMSF}}$ is treeable and $\mathbf P$-almost every orbit is a tree with one end. It follows that $\mathcal R_{\mathsf{WMSF}}$ is hyperfinite. We prove the following elementary fact (see \cite[Proposition 3.5]{lps}).

\begin{prop}\label{ends}
Let $\mathcal G = \Cay(\Gamma, S)$. Assume that $\mathsf{WMSF} \neq \mathsf{FMSF}$. Then $\mathcal R_{\mathsf{FMSF}}$ is not hyperfinite.
\end{prop}

\begin{proof}
Assume that $\mathcal R_{\mathsf{FMSF}}$ is hyperfinite. Using \cite[Proposition ${\rm III}$.3]{cout}, we get $1 \leq \cost(\mathcal R_{\mathsf{WMSF}}) \leq \cost(\mathcal R_{\mathsf{FMSF}}) = 1$ so that $\mathcal R_{\mathsf{WMSF}} = \mathcal R_{\mathsf{FMSF}}$. For $\omega = \mathsf{WMSF}(x)$ or $\mathsf{FMSF}(x)$, denote by $\mathsf{T}(\omega; g)$ the tree (cluster) containing the vertex $g \in \Gamma$. Therefore, $\mathbf P$-a.s.\ $\mathsf T(\mathsf{WMSF}(x); 1_\Gamma) = \mathsf T(\mathsf{FMSF}(x); 1_\Gamma)$. By $\Gamma$-invariance, we get that $\mathbf P$-a.s.\ for all $g \in \Gamma$, $\mathsf T(\mathsf{WMSF}(x); g) = \mathsf T(\mathsf{FMSF}(x); g)$ and thus $\mathsf{WMSF} = \mathsf{FMSF}$.
\end{proof}

Tim\'ar \cite{timar} proved that if $\mathsf{WMSF} \neq \mathsf{FMSF}$, then $\mathcal R_{\mathsf {FMSF}}$ is in fact nowhere hyperfinite, that is, the restriction of $\mathcal R_{\mathsf {FMSF}}$ to any non-null measurable subset is not hyperfinite. We now present the proof of the result of Gaboriau and Lyons \cite{gaboriau-lyons}. We will use a result of Chifan and Ioana \cite[Theorem 1]{CI08}, the proof of which is postponed until Section \ref{subequivalence-bernoulli}.

\begin{theo}[Measurable subgroup]\label{measurable}
For any nonamenable group $\Gamma$ there exists a free ergodic pmp action $\F_2 \curvearrowright ([0, 1]^\Gamma, \Leb^\Gamma)$ such that 
$$\mathcal R(\F_2 \curvearrowright [0, 1]^\Gamma) \subset \mathcal R(\Gamma \curvearrowright [0, 1]^\Gamma).$$
\end{theo}

\begin{proof}
Let $\Gamma$ be a nonamenable group. Since the union of an increasing sequence of amenable groups is still amenable, $\Gamma$ contains a nonamenable finitely generated subgroup. Thus, up to taking such a subgroup, we may assume that $\Gamma$ is finitely generated. The proof is in two steps.

{\bf Step 1.} There exists a subequivalence relation $\mathcal R \subset \mathcal R(\Gamma \curvearrowright [0, 1]^\Gamma)$ which is ergodic treeable and non-hyperfinite.

Let $S$ be a finite generating family such that the Cayley graph $\mathcal G = \Cay(\Gamma, S)$ satisfies $p_c(\mathcal G) < p_u(\mathcal G)$ (see Corollary \ref{pak-smirnova}). As usual, denote the graph $\mathcal G = (\mathsf V, \mathsf E)$. Recall that the pmp actions $\Gamma \curvearrowright [0, 1]^\Gamma$ and $\Gamma \curvearrowright [0, 1]^{\mathsf E}$ are conjugate. By Propositions \ref{forests} and \ref{ends}, we know that $\mathcal R_{\mathsf{FMSF}}$ is not hyperfinite. Apply now Theorem \ref{subequivalence} to $\mathcal R_{\mathsf{FMSF}}$ that we regard as a subequivalence relation of $\mathcal R(\Gamma \curvearrowright [0, 1]^\Gamma)$. Then there exists a non-null measurable subset $X \subset [0, 1]^\Gamma$ such that $\mathcal R_{\mathsf{FMSF}} |X$ is ergodic treeable and non-hyperfinite. In order to extend $\mathcal R_{\mathsf{FMSF}} |X$ to $[0, 1]^\Gamma$, choose an enumeration $\{g_i : i \in \N\}$ of $\Gamma$. For every $x \in [0, 1]^\Gamma \setminus X$, let $n_x$ be the least integer $j \in \N$ such that $g_{j} x \in X$. Let $\mathcal R$ be the smallest equivalence relation containing $\mathcal R_{\mathsf{FMSF}} |X$ and $(x, g_{n_x} x)$, for $x \in [0, 1]^\Gamma \setminus X$. We get that $\mathcal R$ is ergodic treeable and non-hyperfinite.

{\bf Step 2.} There exists a subequivalence relation $\mathcal S \subset \mathcal R(\Gamma \curvearrowright [0, 1]^\Gamma)$ which is induced by a free ergodic pmp action $\F_2 \curvearrowright [0, 1]^\Gamma$.

By \cite[Th\'eor\`eme ${\rm IV}$.1]{cout}, we have that $\mathcal R$ has cost greater than $1$. 
Next, we need the following result due to Hjorth \cite{hjorth-cost} (see also the proof of \cite[Theorem 28.3]{kechris-miller}).

\begin{lemm}\label{hjorth}
Any ergodic treeable pmp equivalence relation $\mathcal R$ such that $\cost(\mathcal R) \geq 2$ contains a subequivalence relation induced by a free pmp action of $\F_2 = \langle a, b\rangle$ such that the generator $a$ acts ergodically.
\end{lemm}

Using the induction formula \cite[Proposition ${\rm II}.6$]{cout}, let $U \subset [0, 1]^\Gamma$ be a Borel measurable subset such that $\cost(\mathcal R | U) \geq 2$. By Lemma \ref{hjorth}, $\mathcal R | U$ contains a subequivalence relation $\mathcal T = \mathcal R(\F_2 \curvearrowright U)$ induced by a free pmp action of $\F_2 = \langle a, b\rangle$ such that the generator $a$ acts ergodically. By considering a subgroup of $\F_2$ of the form $\langle b^kab^{k}: 1 \leq k \leq n\rangle$, for some large $n \in \N$, one gets an ergodic treeable subequivalence relation of $\mathcal R | U$ with large cost so that when extended to the whole space (by using partial Borel isomorphisms of $\mathcal R$), it gets cost $\geq 2$ by \cite[Proposition ${\rm II}.6$]{cout}. Another application of Lemma \ref{hjorth} finishes the proof of Step 2.
\end{proof}

\section{Finite von Neumann algebras}\label{vonNeumann}

We review a few concepts involving finite von Neumann algebras. Further information on this topic may be found in the book \cite{BO} by Brown and Ozawa. 

A {\em von Neumann algebra} $M$ is a unital $\ast$-subalgebra of $\mathbf B(\ell^2)$ which is closed for the strong operator topology. We only deal with {\em tracial} or {\em finite} von Neumann algebras, that is, $M$ is always assumed to carry a faithful normal state $\tau : M \to \C$ which moreover satisfies the {\em trace} identity: $\tau(x y) = \tau(y x)$, for all $x, y \in M$. We denote by $\|x\|_2 = \tau(x^*x)^{1/2}$ the corresponding Hilbert norm and $L^2(M)$ the $L^2$-completion of $M$ with respect to $\|\cdot\|_2$. The uniform norm is denoted by $\|\cdot\|_\infty$. We regard $x \in M$ both as an element of $L^2(M)$ and as a bounded (left multiplication) operator on $L^2(M)$. We will often use the following inequality:
$$\|x \xi y\|_2 \leq \|x\|_\infty \|y\|_\infty \|\xi\|_2, \forall x, y \in M, \forall \xi \in L^2(M).$$
The group of unitaries of $M$ is denoted by $\mathcal U(M)$, the center $M' \cap M$ is $\mathcal Z(M)$ and the unit ball with respect to the uniform norm is $(M)_1$. An infinite dimensional finite von Neumann algebra with trivial center is called a ${\rm II_1}$ {\em factor}.

The main class of examples of finite von Neumann algebras arises from the {\em group measure space construction} of Murray and von Neumann \cite{mvn1}. Let $\Gamma \curvearrowright (X, \mu)$ be a free pmp action of a countable infinite group $\Gamma$ on a nonatomic standard probability space. We regard $F \in L^\infty(X)$ as a bounded operator on $\ell^2(\Gamma) \otimes L^2(X)$ by identifying $F$ with $1 \otimes F \in \mathbf B(\ell^2(\Gamma) \otimes L^2(X))$. The action $\Gamma \curvearrowright X$ induces a unitary representation $\sigma : \Gamma \to \mathcal U(L^2(X))$ defined by $\sigma_g(\xi)(x) = \xi(g^{-1}x)$, for all $\xi \in L^2(X)$. Let $\lambda : \Gamma \to \mathcal U(\ell^2(\Gamma))$ be the left regular representation. The unitaries $u_g = \lambda_g \otimes \sigma_g$ satisfy the following covariance relation: $u_g \xi u_g^* = \sigma_g(\xi)$, for all $\xi \in L^2(X)$, $g \in \Gamma$. By Fell's absorption principle, the unitary representation $(u_g)_{g \in \Gamma}$ is unitarily equivalent to a multiple of $(\lambda_g)_{g \in \Gamma}$. The {\em crossed product} von Neumann algebra $L^\infty(X) \rtimes \Gamma$ is defined by 
$$L^\infty(X) \rtimes \Gamma := \left \{ \sum_{\finite} \xi_g u_g: \xi_g \in L^\infty(X)\right \}'' \subset \mathbf B(\ell^2(\Gamma) \otimes L^2(X)).$$
The von Neumann algebra $M := L^\infty(X) \rtimes \Gamma$ contains a copy of $L^\infty(X)$ as well as a copy of the group von Neumann algebra $L(\Gamma)$. Moreover $M$ is endowed with a trace $\tau$ given by $\tau(a) = \langle a(\delta_e \otimes \mathbf 1_X), \delta_e \otimes \mathbf 1_X\rangle$. The subalgebra $A := L^\infty(X) \subset M$ is called a {\em Cartan subalgebra}.\footnote{A Cartan subalgebra $A \subset M$ is a maximal abelian $\ast$-subalgebra whose normalizer $\mathcal N_M(A) = \{u \in \mathcal U(M) : u A u^* = A\}$ generates $M$ as a von Neumann algebra.} The von Neumann algebra $M$ is a ${\rm II_1}$ factor if and only if the action $\Gamma \curvearrowright X$ is ergodic.
More generally, one can define the von Neumann algebra $L(\mathcal R)$ of a pmp equivalence relation $\mathcal R$ on $(X, \mu)$ (see \cite{feldman-moore2}). Note that $L^\infty(X) \subset L(\mathcal R)$ is still a Cartan subalgebra. When $\mathcal R$ is a type ${\rm II_1}$ equivalence relation, $\mathcal R$ is ergodic if and only if $L(\mathcal R)$ is a ${\rm II_1}$ factor. For a free pmp action $\Gamma \curvearrowright (X, \mu)$, the von Neumann algebras $L^\infty(X) \rtimes \Gamma$ and $L(\mathcal R(\Gamma \curvearrowright X))$ are $\ast$-isomorphic.

Given finite von Neumann algebras $M$ and $N$, an $M$-$N$-{\em bimodule} $\vphantom{}_M\mathcal H_N$ is a Hilbert space endowed with two commuting normal $\ast$-representations $\pi_M : M \to \mathbf B(\mathcal H)$ and $\pi_{N^{\op}} : N^{\op} \to \mathbf B(\mathcal H)$. We simply denote $x \xi y = \pi_M(x) \pi_{N^{\op}}(y)\xi$, for all $x \in M$, $y \in N$, $\xi \in \mathcal H$. The bimodule $\vphantom{}_ML^2(M)_M$ is the {\em trivial bimodule} and $\vphantom{}_{M \otimes 1} L^2(M \overline\otimes M)_{1 \otimes M}$ is the {\em coarse bimodule}. Given two $M$-$N$-bimodules $\mathcal H$ and $\mathcal K$, we say that $\mathcal H$ is {\em weakly contained} in $\mathcal K$ and write $\mathcal H \subset_{\weak} \mathcal K$, if for all $\xi, \eta \in \mathcal H$ and all finite subsets $F \subset M$, $G \subset N$, there exist two sequences $\xi_n, \eta_n$ in finite direct sums of $\mathcal K$ such that
$$\langle x \xi y, \eta \rangle = \lim_n \langle x \xi_n y, \eta_n \rangle, \forall x \in F, \forall y \in G.$$

Given an inclusion $B \subset M$ of finite von Neumann algebras, denote by $E_B : M \to B$ the unique trace-preserving normal {\em conditional expectation}. If we moreover denote by $e_B : L^2(M) \to L^2(B)$ the orthogonal projection, we have $e_B x e_B = E_B(x) e_B$, for all $x \in M$. The basic construction $\langle M, e_B\rangle$ is the von Neumann subalgebra of $\mathbf B(L^2(M))$ generated by $M$ and $e_B$. It is endowed with a faithful normal semifinite trace $\Tr$ given by $\Tr(x e_B y) = \tau(x y)$, for all $x, y \in M$. The $M$-$M$-bimodule $L^2(\langle M, e_B\rangle)$ is {\em mixing relative to} $B$ in the following sense: whenever $u_n \in \mathcal U(M)$ is a sequence of unitaries such that $\lim_n \|E_B(x^* u_n y)\|_2 = 0$, for all $x, y \in M$, then for every $\xi, \eta \in L^2(\langle M, e_B\rangle)$, we have
$$\lim_n \sup_{y \in (M)_1} \left| \langle u_n \xi y, \eta \rangle \right | = \lim_n \sup_{x \in (M)_1} \left| \langle x \xi u_n, \eta \rangle \right | = 0.$$

Recall that $M$ is {\em hyperfinite} if there exists an increasing sequence of unital finite dimensional $\ast$-subalgebras $Q_n \subset M$ such that $M$ is the weak closure of $\bigcup_{n} Q_n$. When $\mathcal R$ is a pmp equivalence relation, $\mathcal R$ is hyperfinite if and only if $L(\mathcal R)$ is hyperfinite \cite{cfw}. In their seminal work \cite{mvn}, Murray and von Neumann showed the uniqueness of the hyperfinite ${\rm II_1}$ factor. We say that $M$ is {\em amenable} if 
$$\vphantom{}_ML^2(M)_M \subset_{\weak} \vphantom{}_{M \otimes 1}L^2(M \overline\otimes M)_{1 \otimes M}.$$
Any hyperfinite von Neumann algebra is amenable. By Connes' groundbreaking work \cite{connes76}, any amenable von Neumann algebra is hyperfinite. Therefore, there is a unique amenable ${\rm II_1}$ factor.

Recall at last Popa's intertwining-by-bimodules technique. Popa discovered \cite{{popamal1}, {popa2001}} a very powerful technique to unitarily conjugate subalgebras in an ambient von Neumann algebra. Let $A, B \subset M$ be subalgebras of a finite von Neumann algebra. The following are equivalent (see \cite[Theorem 2.1]{popamal1}, \cite[Theorem A.1]{popa2001} and also \cite[Theorem C.3]{vaesbern}).
\begin{itemize}
\item There exist projections $p \in A$, $q \in B$, a nonzero partial isometry $v \in pMq$ and a $\ast$-homomorphism $\varphi : pAp \to qBq$ such that $x v = v \varphi(x)$, for all $x \in pAp$.
\item There is no sequence of unitaries $u_n \in \mathcal U(A)$ such that 
$$\lim_n \|E_B(x u_n y)\|_2 = 0, \forall x, y \in M.$$
\end{itemize}
If one of the two conditions holds, we say that $A$ {\em embeds into} $B$ {\em inside} $M$ and write $A \preceq_M B$. By definition, $A$ is {\em diffuse} if $A \npreceq_A \C$, that is, if $A$ has no nonzero minimal projection.

\section{Subequivalence relations of Bernoulli actions}\label{subequivalence-bernoulli}

As we have seen before, given a Cayley graph $\mathcal G = \Cay(\Gamma, S)$, a $\Gamma$-equivariant map $\pi : [0, 1]^\mathsf{E} \to \{0, 1\}^{\mathsf E}$ gives rise to a percolation $\pi_\ast \mathbf P$ on $\mathcal G$ and hence to a subequivalence relation $\mathcal R_\pi^{\cl}$ of the equivalence relation $\mathcal R(\Gamma \curvearrowright [0, 1]^{\mathsf E})$ induced by the Bernoulli action. The aim of this section is to present a global dichotomy result for subequivalence relations of $\mathcal R(\Gamma \curvearrowright [0, 1]^{\mathsf E})$, obtained by Chifan and Ioana \cite[Theorem 1]{CI08}.

\begin{theo}[Dichotomy for subequivalence relations]\label{subequivalence}
Let $\Gamma$ be any infinite countable discrete group. Let $\mathcal R \subset \mathcal R(\Gamma \curvearrowright [0, 1]^\Gamma)$ be any subequivalence relation of the pmp equivalence relation induced by the Bernoulli action. Then there exists a measurable partition $\{X_n : n \in \N \}$ of $[0, 1]^\Gamma$ into $\mathcal R$-invariant subsets such that 
\begin{itemize}
\item $\mathcal R | X_0$ is hyperfinite.
\item $\mathcal R | X_n$ is strongly ergodic, for all $n \geq 1$.
\end{itemize}
\end{theo}

We give a self-contained proof of this result. We first start by recalling the construction of the {\em support length deformation} for Bernoulli actions due to Ioana \cite{ioana06}. We will be using the following notation throughout this section. 
\begin{itemize}
\item Let $(A_0, \tau)$ be an abelian von Neumann algebra, $A = A_0^\Gamma$ the infinite tensor product indexed by $\Gamma$ and $\Gamma \curvearrowright A$ the corresponding Bernoulli shift. Set $M = A \rtimes \Gamma$.
\item Likewise, let $B_0 = A_0 \ast L(\Z)$ be the free product with respect to the natural traces, $B = B_0^\Gamma$ and $\sigma : \Gamma \curvearrowright B$ the corresponding Bernoulli shift. Set $\widetilde M = B \rtimes \Gamma$.
\end{itemize}
Observe that $M \subset \widetilde M$ and denote by $E_M : \widetilde M \to M$ the unique trace-preserving normal conditional expectation.
Following \cite{ioana06}, denote by $v \in L(\Z)$ the canonical generating Haar unitary and take the selfadjoint element $h \in L(\Z)$ with spectrum $[-\pi, \pi]$ such that $v = \exp(i h)$. Denote by $\theta_t^0 \in \Aut(B_0)$ the inner automorphism given by $\theta_t^0 = \Ad(\exp(i t h))$ and let $\theta_t = \otimes_{g \in \Gamma} \theta_t^0 \in \Aut(B)$. Since $(\theta_t)$ commutes with the Bernoulli action, we can extend $(\theta_t)$ to $\widetilde M$ by letting $\theta_t(u_g) = u_g$. We get that $(\theta_t)_{t \in \R}$ is a one-parameter group of automorphisms of $\widetilde M$ such that $\lim_{t \to 0} \|x - \theta_t(x)\|_2 = 0$, for all $x \in M$. Denote by $\beta_0 \in \Aut(B_0)$ the automorphism given by $\beta_0(a) = a$, for all $a \in A_0$ and $\beta_0(v) = v^*$. Define $\beta = \otimes_{g \in \Gamma} \beta_0$ and extend $\beta$ to $\widetilde M$ by acting trivially on $L(\Gamma)$. By construction, $\beta | M = \Id_M$, $\beta^2 = \Id_{\widetilde M}$  and $\beta \circ \theta_t = \theta_{-t} \circ \beta$, for all $t \in \R$. 

For $0 < \rho < 1$,  define the {\em support length deformation} ${\rm m}_\rho : M \to M$ by
$${\rm m}_\rho(a u_g) = \rho^n a u_g, \forall g \in \Gamma, \forall a \in (A_0 \ominus \C1)^J, J \subset \Gamma, |J| = n.$$
Let $\rho_t = |\sin(\pi t)|^2/|\pi t|^2$. One checks that $(E_M \circ \theta_t)(x) = {\rm m}_{\rho_t}(x)$, for all $x \in M$. In particular, $({\rm m}_\rho)$ is a family of trace-preserving unital completely positive maps for which $\theta_t : M \to \widetilde M$ is a {\em dilation}. In this respect, the support length deformation $({\rm m}_\rho)$ is a variant of the malleable deformation discovered by Popa in \cite{popamal1}. Popa used his malleable deformation together with his intertwining techniques to prove various striking rigidity results for Bernoulli actions (see for instance \cite{{popasup}, {popamal1}} and Vaes' Bourbaki seminar \cite{vaesbern} on this topic.)

Spectral gap rigidity was discovered by Popa \cite{{popasup}, {popasolid}}. It was a completely new type of rigidity where the usual (relative) property (T) assumption in many (orbit and W$^*$)-rigidity results could be dropped. Using this technique, Popa \cite{popasup} proved, among other results, that for any nonamenable product of infinite groups $\Gamma = \Gamma_1 \times \Gamma_2$, the plain Bernoulli action $\Gamma \curvearrowright [0, 1]^\Gamma$ is $\mathcal U_{\fin}$-cocycle superrigid.\footnote{$\mathcal U_{\fin}$ is the class of groups which embed into the unitary group of a $\|\cdot\|_2$-separable ${\rm II_1}$ factor.}

The following variant of spectral gap property is due to Chifan and Ioana (see \cite[Lemma 5]{CI08}).

\begin{prop}[Spectral gap]\label{spectralgap}
As $M$-$M$-bimodules, we have
\begin{equation}\label{weak}
\vphantom{}_M(L^2(\widetilde M) \ominus L^2(M))_M \subset_{\weak} \vphantom{}_{M \otimes 1}L^2(M \overline \otimes M)_{1 \otimes M}.
\end{equation}
\end{prop}

\begin{proof}
We start by proving the following.
\begin{claim}
There is a countable set $\{(\Gamma_i, \Delta_i) : i \in \mathcal I\}$, where $\Gamma_i < \Gamma$ is a finite subgroup and $\Delta_i \subset \Gamma$ is a non-empty set which is invariant under left multiplication by $\Gamma_i$ such that with $A_i = A_0^{\Gamma \setminus \Delta_i} \rtimes \Gamma_i$, we have an isomorphism of $M$-$M$-bimodules
\begin{equation}\label{bimodule}
L^2(\widetilde M) \ominus L^2(M) \cong \bigoplus_{i \in \mathcal I} L^2 \left ( \left\langle M, e_{A_i} \right\rangle \right).
\end{equation}
\end{claim}

To prove the claim, let $\mathcal A_0 \subset A_0 \ominus \C$ be an orthonormal basis of $L^2(A_0) \ominus \C$ and denote by $v$ the Haar unitary generating $L(\mathbf Z)$. Recall that $B_0 = A_0 \ast L(\Z)$. Define the subset $\mathcal B_0 := \{v^{n_1} a_1 \cdots v^{n_k} a_k v^{n_{k +  1}} : k \geq 0, n_1, \dots, n_{k +  1} \in \mathbf Z - \{0\}, a_i \in \mathcal A_0\}$.
By construction, we have a decomposition
$$L^2(B_0) \ominus L^2(A_0) = \bigoplus_{b \in \mathcal B_0} \overline{A_0 b A_0}$$
into pairwise orthogonal $A_0$-$A_0$-subbimodules. Define the countable set 
$$\mathcal I = \left\{b_{\mathcal F} = \bigotimes_{g \in \mathcal F} b_g : \mathcal \emptyset \neq \mathcal F \subset \Gamma \mbox{ finite subset}, b_g \in \mathcal B_0 \mbox { for all } g \in \mathcal F \right\}.$$ We have a decomposition
\begin{equation}\label{prebimodule}
L^2(\widetilde M) \ominus L^2(M) = \bigoplus_{b \in \mathcal I} \overline{M b M}
\end{equation}
into pairwise orthogonal $M$-$M$-subbimodules. For $b \in \mathcal I$, define the finite subgroup $\Gamma_b = \{g \in \Gamma : g\mathcal F = \mathcal F \mbox{ and } \sigma_g(b) = b\}$. Let $A_b = A_0^{\Gamma \setminus \mathcal F} \rtimes \Gamma_b$. One checks that the map $x e_{A_b} y \mapsto x b y$ defines an $M$-$M$-bimodule isomorphism
\begin{equation}\label{isomorphism}
L^2(\langle M, e_{A_b}\rangle) \to \overline{M b M}.
\end{equation}
The Claim follows now from $(\ref{prebimodule})$ and $(\ref{isomorphism})$. Finally, since $A_i$ is amenable, the isomorphism $(\ref{bimodule})$ together with \cite[Lemma 1.7]{anan95} yield $(\ref{bimodule})$
\end{proof}

If $P \subset M$ has no amenable direct summand, then for every $\varepsilon > 0$, there exists $\delta > 0$ and $\mathcal V \subset \mathcal U(P)$ finite subset such that for every $x \in (\widetilde M)_1$,
\begin{equation}\label{gap}
(\|u x - x u\|_2 \leq \delta, \forall u \in \mathcal V) \Longrightarrow \|x - E_M(x)\|_2 \leq \varepsilon.
\end{equation}

Indeed, assume that $(\ref{gap})$ does not hold. Then one can find a sequence $x_n \in (M)_1$, such that $x_n \in L^2(\widetilde M) \ominus L^2(M)$, $\| x_n \|_2 = 1$ and $\lim_n \|y x_n - x_n y\|_2 = 0$, for all $y \in P$. Up to passing to a subsequence we may assume that $b_n = x_n x_n^*$ converges weakly to $b \in (P' \cap M)_1$. Observe that $\tau(b) = 1$. Let $c \in \mathcal Z(P)_+$ so that $p = E_P(b)^{1/2} c \in \mathcal Z(P)$ is a nonzero projection. From $(\ref{weak})$, we get that, as $Pp$-$Pp$-bimodules,
\begin{equation}\label{bimodule2}
\vphantom{}_{Pp}(L^2(\widetilde M) \ominus L^2(M))_{Pp} \subset_{\weak} \vphantom{}_{Pp \otimes 1}L^2(Pp \overline\otimes Pp)_{1 \otimes Pp}.
\end{equation}
Define $\xi_n := c x_n$. For all $y \in P$, we have $\lim_n \|y \xi_n - \xi_n y\|_2 = 0$ and 
$$\lim_n \langle y \xi_n, \xi_n\rangle =  \lim_n \tau(y c x_n x_n^* c) = \lim_n\tau(y c b c) =  \tau(yp),$$
whence 
\begin{equation}\label{bimodule3}
\vphantom{}_{Pp}L^2(Pp)_{Pp} \subset_{\weak} \vphantom{}_{Pp}(L^2(\widetilde M) \ominus L^2(M))_{Pp}.
\end{equation}
Together with $(\ref{bimodule2})$ and $(\ref{bimodule3})$, we finally obtain that $Pp$ is amenable.

The next result due to Chifan and Ioana (see \cite[Theorem 2]{CI08}) is the key to proving the global dichotomy result for subequivalence relations. 

\begin{theo}\label{commutant}
Let $Q \subset A$ be a diffuse von Neumann subalgebra. Then $Q' \cap M$ is amenable.
\end{theo}

We point out that this result was earlier obtained by Ozawa \cite[Theorem 4.7]{ozawa2004} for all {\em exact} groups $\Gamma$ using C$^*$-algebraic techniques. Chifan and Ioana's proof that we present here relies on a theory developed by Popa over the last decade known today as deformation vs.\ rigidity. We refer to \cite{{popa-icm}, {vaes-icm}} for further information on this topic.

\begin{proof}[Proof of Theorem $\ref{commutant}$]
The proof is reminiscent of the one of \cite[Theorem 4.1]{popamal1} (see also \cite[Lemma 6.1]{vaesbern}). We prove the result by contradiction following the lines of the proof of \cite[Theorem 4.2]{IPV}. We may assume that $Q \subset A$ is diffuse and $Q' \cap M$ has no amenable direct summand. We will be using the following terminology. Given subalgebras $Q_1, Q_2 \subset \widetilde M$, an element $x \in \widetilde M$ is said to be $Q_1$-$Q_2$-finite inside $\widetilde M$ if there exist elements $x_1, \dots, x_m, y_1, \dots, y_n \in \widetilde M$ such that
\begin{equation}\label{finite}
x Q_2 \subset \sum_{i = 1}^m Q_1 x_i \mbox{ and } Q_1 x \subset \sum_{j = 1}^n y_j Q_2.
\end{equation}

\begin{newstep}
There exist $t = 1/2^n$ and a nonzero element $v \in \widetilde M$ which is $Q$-$\theta_t(Q)$-finite.
\end{newstep}

Let $\varepsilon  = 1/2$. Proposition \ref{spectralgap} yields $\delta > 0$ and a finite subset $\mathcal V \subset \mathcal U(Q' \cap M)$ for which $(\ref{gap})$ holds. Let $s$ small enough so that $\|b - \theta_s(b)\|_2 \leq \delta/2$, for all $b \in \mathcal V$. For all $u \in \mathcal U(Q)$,
\begin{eqnarray*}
\|b \theta_s(u) - \theta_s(u) b\|_2 & = & \|(b - \theta_s(b)) \theta_s(u) - \theta_s(u) (b - \theta_s(b))\|_2 \\
& \leq & 2 \|\theta_s(u)\|_\infty \|b - \theta_s(b)\|_2 \leq \delta.
\end{eqnarray*}
Using Proposition \ref{spectralgap}, we get $\|\theta_s(u) - E_M(\theta_s(u))\|_2 \leq 1/2$, for all $u \in \mathcal U(Q)$. Let $\rho = \rho_s^2$, so that ${\rm m}_\rho = {\rm m}_{\rho_s}^2$. For all $u \in \mathcal U(Q)$, we have
$$1 - \tau(u^*{\rm m}_{\rho}(u)) = 1 - \|{\rm m}_{\rho_s}(u)\|_2^2 = \|\theta_s(u) - E_M(\theta_s(u))\|_2^2 \leq 1/4.$$
Then $\tau(u^* \theta_s(u)) = \tau(u^* {\rm m}_\rho(u)) \geq 3/4$, for all $u \in \mathcal U(Q)$. Since $t \mapsto \tau(u^* \theta_t(u))$ is decreasing, we can take $t = 1/2^n$ such that $\tau(u^*\theta_t(u)) \geq 3/4$, for all $u \in \mathcal U(Q)$. 
Let $v$ be the unique element of minimal $\|\cdot\|_2$-norm in the weak closure of the 
convex hull of $\{u^* \theta_t(u) : u \in \mathcal U(Q)\}$. We get $\tau(v) \geq 3/4$ and $u v =  v \theta_t(u)$, for all $u \in \mathcal U(Q)$ (by uniqueness). In particular, $v \in \widetilde M$ is a nonzero $Q$-$\theta_t(Q)$-finite element.

\begin{newstep}
There exists a nonzero element $a \in \widetilde M$ which is $Q$-$\theta_1(Q)$-finite.
\end{newstep}

To prove Step $2$, it suffices to show the following statement: if there exists a nonzero element $v$ which is $Q$-$\theta_t(Q)$-finite, then there exists a nonzero element $w$ which is $Q$-$\theta_{2t}(Q)$-finite. Indeed, since $t = 1/2^n$, we can then go until $t = 1$. Denote by $\QN_M(Q)$ the set of all $Q$-$Q$-finite elements inside $M$ ($\QN_M(Q)$ is also called the {\em quasi-normalizer of} $Q$ {\em inside} $M$ \cite{popa2001}). Let $P := \QN_M(Q)'' \subset M$. Observe that for all $d \in \QN_M(Q)$, the element $\theta_t(\beta(v^*) d v)$ is $Q$-$\theta_{2t}(Q)$-finite. Indeed, let $d \in \QN_M(Q)$ which satisfies $(\ref{finite})$ for $Q_1 = Q_2 = Q$. Then we get
\begin{eqnarray*}
\theta_t(\beta(v^*) d v) \theta_{2t}(Q) & = & \theta_t(\beta(v^*) d Q v) \subset \sum_i \theta_t(\beta(v^*) Q x_i v) = \sum_i Q \theta_t(\beta(v^*) x_i v) \\
Q\theta_t(\beta(v^*) d v) & = & \theta_t(\beta(v^*) Q d v) \subset \sum_j \theta_t(\beta(v^*)y_j Q v) = \sum_j \theta_t(\beta(v^*)y_j v) \theta_{2t}(Q).
\end{eqnarray*}
Hence we have to prove that there exists $d \in \QN_M(Q)$ such that $\beta(v^*) d v \neq 0$. By contradiction, assume that this is not the case. Denote by $q \in \widetilde M$ the projection onto the closed linear span of $\{\range(dv) : d \in \QN_M(Q)\}$. We have $\beta(v^*)q = 0$ and $q \in P' \cap \widetilde M$.

We use now again the $M$-$M$-bimodule isomorphism $(\ref{bimodule})$. Since $Q' \cap M \subset P$, it follows that $P$ has no amenable direct summand and thus $P \npreceq_M A_i$, for all $i \in \mathcal I$. Therefore there exists a sequence of unitaries $u_n \in \mathcal U(P)$ such that $\lim_n \|E_{A_i}(x^* u_n y)\|_2 = 0$, for all $x, y \in M$, $i \in I$. Let $x \in P' \cap \widetilde M$. Set $\eta := x - E_M(x)$. Observe that $\eta \in P' \cap \widetilde M$ and $\eta \perp L^2(M)$. Write $\eta = \oplus_{i \in \mathcal I} \eta_i$, with $\eta_i \in L^2(\langle M, e_{A_i}\rangle)$. Since the $M$-$M$-bimodule $L^2(\langle M, e_{A_i}\rangle)$ is mixing relative to $A_i$, we have $\lim_n \langle u_n \eta_i u_n^*, \eta_i \rangle = 0$, for all $i \in \mathcal I$ and so $\lim_n \langle u_n \eta u_n^*, \eta \rangle = 0$. Since $\eta \in P' \cap \widetilde M$, we have $\|\eta\|^2_2 = \lim_n \langle u_n \eta u_n^*, \eta \rangle = 0$. Therefore $P' \cap \widetilde M = P' \cap M$. In particular, we get $q \in M$, so that $\beta(v^* q) = \beta(v^*)q = 0$. Hence $v = 0$, which is a contradiction. 

Observe that $\overline{M a \theta_1(Q)}$ is a nonzero $M$-$\theta_1(Q)$-subbimodule of $L^2(\widetilde M)$ which is finitely generated as left $M$-module, whence we get $\theta_1(Q) \preceq_{\widetilde M} M$. We use the following notation: for every nonempty finite subset $\mathcal F \subset \Gamma$, let $\Stab(\mathcal F) = \{g \in \Gamma : g\mathcal F = \mathcal F\}$ and $M(\mathcal F) := A_0^\mathcal{F} \rtimes \Stab(\mathcal F)$. By convention, set $M(\emptyset) := L(\Gamma)$.

\begin{newstep}
There exists a finite subset $\mathcal F \subset \Gamma$ such that $Q \preceq_M M(\mathcal F)$.
\end{newstep}

We prove Step $3$ by contradiction and assume that for all finite subset $\mathcal F \subset \Gamma$, we have $Q \npreceq_M M(\mathcal F)$. Let $v_n \in \mathcal U(Q)$ be a sequence of unitaries such that $\lim_n \| E_{M(\mathcal F)}(x^* v_n y) \|_2 = 0$, for all $x, y \in M$, $\mathcal F \subset \Gamma$. We upgrade this by showing the following:
\begin{equation}\label{convergence}
\lim_n \|E_M(x^* \theta_1(v_n) y)\|_2 = 0, \forall x, y \in \widetilde M.
\end{equation}
This will clearly contradict Step 2. Let $\mathcal F, \mathcal G \subset \Gamma$ be finite (possibly empty) subsets. Define $x = \bigotimes_{g \in \mathcal F} x_g \otimes \bigotimes_{g \in \Gamma \setminus \mathcal F} 1$ and $y = \bigotimes_{h \in \mathcal G} y_h \otimes \bigotimes_{h \in \Gamma \setminus \mathcal G} 1$, where $x_g, y_h \in B_0 \ominus \theta_1(A_0) A_0$. Observe that it suffices to prove $(\ref{convergence})$ for such $x$ and $y$ since the linear span of all $\theta_1(A) y M$ for $y$ of the above form is a $\|\cdot\|_2$-dense subspace of $\widetilde M$.

Write $v_n = \sum_{g \in \Gamma} (v_n)^g u_g$ for the Fourier expansion of $v_n$ in $M$, where $(v_n)^g \in A$. We have $E_M(x^* \theta_1(v_n) y) = \sum_{g \in \Gamma} E_A\left(x^* \theta_1((v_n)^g)\sigma_g(y)\right) u_g$. If $g\mathcal G \neq \mathcal F$, then $E_A\left( x^* \theta_1((v_n)^g)\sigma_g(y) \right) = 0$. If $g\mathcal G = \mathcal F$, then
$$E_A\left( x^* \theta_1((v_n)^g)\sigma_g(y) \right) = E_A\left( x^* \theta_1\left( E_{A_0^\mathcal{F}}((v_n)^g)\right) \sigma_g(y) \right).$$
Take now finitely many $g_1, \dots, g_k \in \Gamma$ such that $g_i\mathcal G = \mathcal F$ and such that  $\{g \in \Gamma : g\mathcal G = \mathcal F\}$ is the disjoint union of $(\Stab \mathcal F)g_1, \dots, (\Stab\mathcal F)g_k$. Set $w_n = \sum_{i = 1}^k E_{M(\mathcal F)} (v_n u_{g_i}^*) u_{g_i}$. We have proven $E_M(x^* \theta_1(v_n) y) = E_M(x^* \theta_1(w_n) y)$.
Since by assumption $\lim_n \|w_n\|_2 = 0$, we get $(\ref{convergence})$.

\begin{newstep}
We derive a contradiction.
\end{newstep}

From Step 3, there exists a finite subset $\mathcal F \subset \Gamma$ such that $Q \preceq_M M(\mathcal F)$. If $\mathcal F = \emptyset$, then $Q \preceq_M L(\Gamma)$. Since $M = A \rtimes \Gamma$, this clearly contradicts the fact that $Q \subset A$ is diffuse. Hence $\mathcal F \neq \emptyset$ and since $\Stab(\mathcal F)$ is finite, we get $Q \preceq_M A_0^\mathcal{F}$. There exist projections $q \in Q$, $r \in A_0^\mathcal{F}$, a nonzero partial isometry $v \in qMr$ and a $\ast$-homomorphism $\varphi : qQq \to r A_0^\mathcal{F}r$ such that $x v = v \varphi(x)$, for all $x \in qQq$. Hence $\varphi(qQq) \subset r A_0^\mathcal F r$ is a diffuse subalgebra. A straightforward computation shows that $\varphi(qQq)' \cap r Mr \subset r(\sum_{g \in \mathcal G} Au_g)r$, where $\mathcal G = \mathcal F \mathcal F^{-1}$. Since $v^* (Q' \cap M) v \subset \varphi(qQq)' \cap rMr$, we get $v^* (Q' \cap M) v \subset r(\sum_{g \in \mathcal G} Au_g)r$. Thus $Q' \cap M \preceq_M A$, which contradicts the fact that $Q' \cap M$ has no amenable direct summand. The proof is complete.
\end{proof}

\begin{proof}[Proof of Theorem $\ref{subequivalence}$]
Let $\mathcal R \subset \mathcal R(\Gamma \curvearrowright [0, 1]^\Gamma)$ be any pmp subequivalence relation. Write $N = L(\mathcal R)$ for the von Neumann algebra of $\mathcal R$. Denote by $z_0 \in \mathcal Z(N)$ the maximal central projection for which $N z_0$ is amenable. We claim that $\mathcal Z(N)(1 - z_0)$ is purely atomic. Assume that this is not the case. Let $q \in \mathcal Z(N)(1 - z_0)$ be a nonzero projection such that $\mathcal Z(N)q$ is diffuse. Set $Q := A(1 - q) \oplus \mathcal Z(N)q \subset A$, which is a diffuse von Neumann subalgebra of $A$. Theorem \ref{commutant} implies that $Q' \cap M$ is amenable and thus $N q$ is amenable, which contradicts the maximality of $z_0$.

Write $\mathcal Z(N)(1 - z_0) = \bigoplus_{n \geq 1} \C z_n$. Denote by $X_n \subset [0, 1]^\Gamma$ the measurable $\mathcal R$-invariant subset corresponding to the central projection $z_n$, that is, $\mathbf 1_{X_n} = z_n$ and $L(\mathcal R | X_n) = N z_n$. We get that $\mathcal R | X_0$ is hyperfinite and $\mathcal R | X_n$ is ergodic and non-hyperfinite, for all $n \geq 1$. In particular, it follows that any subequivalence $\mathcal T \subset \mathcal R(\Gamma \curvearrowright [0, 1]^\Gamma)$ which has a diffuse ergodic decomposition must be hyperfinite. Furthermore, we deduce that $\mathcal R | X_n$ cannot be written as an increasing union of subequivalence relations with a diffuse ergodic decomposition (otherwise $\mathcal R | X_n$ would be hyperfinite). Using Proposition \ref{strong-ergodicity}, we finally obtain that $\mathcal R | X_n$ is strongly ergodic, for all $n \geq 1$.
\end{proof}

\section{Co-induced actions}\label{induction}
Ioana \cite{ioana07} used the co-induction technique \cite{gaboriau-ME} together with a separability argument (see Theorem \ref{separability}) to prove that any nonamenable group $\Gamma$ that contains $\F_2$ has uncountably many non-orbit equivalent actions. First recall the co-induction construction for a subgroup $\Lambda < \Gamma$. Let $\alpha : \Lambda \curvearrowright (Y, \nu)$ be any free pmp action on the nonatomic standard probability space. Fix a section $s : \Gamma/\Lambda \to \Gamma$ such that $s(\Lambda) = 1_\Gamma$. Define the $1$-cocycle $\omega : \Gamma \times \Gamma/\Lambda \to \Lambda$ by $\omega(g, t) = s(gt)^{-1} g s(t)$. The {\em co-induced action} $\sigma = \coInd_\Lambda^\Gamma(\alpha) : \Gamma \curvearrowright (Y^{\Gamma/\Lambda}, \nu^{\Gamma/\Lambda})$ is then defined by $(\sigma_g(y))_t = \alpha(\omega(g, g^{-1} t))(y_{g^{-1}t})$, for all $g \in \Gamma$, $t \in \Gamma/\Lambda$. In order to prove that any nonamenable group has uncountably many non-orbit equivalent actions, we review now Epstein's construction \cite{epstein} of the co-induced action for a measurable subgroup $\Lambda <_{\ME} \Gamma$. 

Let $a : \Lambda \curvearrowright (X, \mu)$ and $b : \Gamma \curvearrowright (X, \mu)$ be free ergodic pmp actions of infinite countable discrete groups $\Lambda$ and $\Gamma$ on the nonatomic standard probability space $(X, \mu)$ such that $\mathcal R(a, \Lambda) \subset \mathcal R(b, \Gamma)$. We will assume that $\mathcal R(a, \Lambda)$ has infinite index in $\mathcal R(b, \Gamma)$, that is, $\mu$-almost every $\mathcal R(b, \Gamma)$-class contains infinitely many $\mathcal R(a, \Lambda)$-classes. Fix {\em choice} functions $(C_n : X \to X)_{n \in \N}$ so that every $C_n : X \to X$ is Borel; $C_0 = \Id_X$; given $x \in X$, $\{C_n(x) : n \in \N\}$ enumerates a tranversal for the $\mathcal R(a, \Lambda)$-classes in the $\mathcal R(b, \Gamma)$-class of $x$; and for all $m \neq n$ and $x \in X$, we have $C_m(x) \neq C_n(x)$. Observe that since $a$ is ergodic, we may assume that the choice functions $C_n$ are one-to-one. 

Denote by $S_\infty$ the full permutation group of $\N$. Let $\mathbf i : \Gamma \times X \to S_\infty$ be the {\em index} cocycle given by the formula
$$\mathbf i (g, x)(k) = n \Longleftrightarrow [C_k(x)]_{\mathcal R(a, \Lambda)} = [C_n(gx)]_{\mathcal R(a, \Lambda)}.$$
Since the action $a : \Lambda \curvearrowright X$ is assumed to be free, we can then define the Borel map $ \ell : \Gamma \times X \to \Lambda^\N$ by the formula
$$\mathbf{\ell}(g, x)_n \cdot C_{\mathbf i (g, x)^{-1}(n)}(x) = C_n(gx).$$
Observe that $S_\infty$ acts on $\Lambda^\N$ by Bernoulli shift: for all $\pi \in S_\infty$ and $(\lambda_n)_{n \in \N} \in \Lambda^\N$, we have $(\pi \cdot \lambda)_n = \lambda_{\pi^{-1}(n)}$. Denote by $S_\infty \ltimes \Lambda^\N$ the corresponding semi-direct product group. We finally define the Borel cocycle $\Omega : \Gamma \times X \to S_\infty \ltimes \Lambda^\N$ by the formula
$$\Omega(g, x) = (\mathbf i(g, x), \ell(g, x)).$$
One checks that $\Omega$ satisfies the $1$-cocycle relation: for $\mu$-almost every $x \in X$, for all $g, h \in \Gamma$, we have $\Omega(gh, x) = \Omega(g, h x) \Omega(h, x)$.

Let now $\alpha : \Lambda \curvearrowright (Y, \nu)$ be any free pmp action on the nonatomic standard probability space. Using the Borel cocycle $\Omega$, we can define the pmp skew-product action $\sigma : \Gamma \curvearrowright (X \times Y^\N, \mu \times \nu^\N)$ by the formula
\begin{eqnarray}\label{skewproduct}
g^\sigma \cdot (x, (y_n)_{n \in \N}) & = & \left(g \cdot x, \Omega(g, x)^{\alpha^\N} \cdot (y_n)_{n \in \N}\right) \\ \nonumber
& = & \left(g \cdot x, \left( n \mapsto (\ell(g, x)_n )^\alpha \cdot y_{\mathbf i (g, x)^{-1}(n)} \right) \right).
\end{eqnarray}
One checks that this action is independent of the choice of $(C_n)_{n \in \N}$, up to conjugation.

\begin{defi}[Co-induced action]
Under the previous assumptions, we say that $\sigma$ is the {\em co-induced action of} $\alpha$ {\em modulo} $(a, b)$ and write 
$$\sigma = \coInd(a, b)_\Lambda^\Gamma(\alpha) : \Gamma \curvearrowright (X \times Y^\N, \mu \times \nu^\N).$$
\end{defi}

We can view $\coInd(a, b)_\Lambda^\Gamma$ as an operation from the space $A(\Lambda , Y, \nu)$ of pmp actions of $\Lambda$ on $(Y, \nu)$ to the space $A(\Gamma, X \times Y^\N, \mu \times \nu^\N)$ (see \cite{kechris}). Observe that when regarding $\Omega : \mathcal R(\Gamma \curvearrowright X) \to S_\infty \ltimes \Lambda^\N$ as a cocycle for the equivalence relation and taking the restriction $\Omega | \mathcal R(\Lambda \curvearrowright X)$, the formula $(\ref{skewproduct})$ also allows to define a skew-product action $\rho : \Lambda \curvearrowright (X \times Y^\N, \mu \times \nu^\N)$ that we will denote by $\rho = \coInd(a, b)_\Lambda^\Lambda(\alpha)$. The action $\rho$ generates a subequivalence relation of the one generated by $\sigma = \coInd(a, b)_\Lambda^\Gamma(\alpha)$, that is, $\mathcal R(\rho, \Lambda) \subset \mathcal R(\sigma, \Gamma)$. Note that 
\begin{itemize}
\item $b$ is a quotient of $\sigma$ with quotient map $(x, (y_n)_{n \in \N}) \mapsto x$.
\item $\alpha$ is a quotient of $\rho$ with quotient map $p_\rho : (x, (y_n)_{n \in \N}) \mapsto y_0$. 
\end{itemize}

In particular, $\rho$ and $\sigma$ are free pmp actions. It turns out that proving ergodicity for the co-induced action $\sigma = \coInd(a, b)_\Lambda^\Gamma(\alpha)$ is more technical and delicate than in the case of a genuine subgroup $\Lambda < \Gamma$. Epstein finds an ergodic measure for the co-induced action $\sigma$ by analyzing the ergodic decomposition of $X$ with respect to the action $b : \Gamma \curvearrowright X$ (see \cite[Lemma 2.6]{epstein}). In \cite{IKT}, Ioana, Kechris and Tsankov circumvent this difficulty by finding necessary and sufficient conditions on the inclusion $\mathcal R(a, \Lambda) \subset \mathcal R(b, \Gamma)$ which ensure that the co-induced action $\sigma$ is mixing, and so ergodic. More precisely, they obtained the following result (see \cite[Theorem 3.3]{IKT}).

\begin{theo}[Mixing co-induced actions]\label{mixing}
Let $a : \Lambda \curvearrowright (X, \mu)$ and $b : \Gamma \curvearrowright (X, \mu)$ be free pmp actions such that $b$ is mixing and $\mathcal R(a, \Lambda) \subset \mathcal R(b, \Gamma)$. Let $N = L^\infty(X) \rtimes_a \Lambda$ and $M = L^\infty(X) \rtimes_b \Gamma$ be the corresponding group measure space von Neumann algebras so that $N \subset M$. Write $(u_g)_{g \in \Gamma}$ for the unitaries in $M$ implementing the action $b$. Denote by $E_N : M \to N$ the trace-preserving normal conditonal expectation. The following are equivalent:
\begin{itemize}
\item $\lim_{g \to \infty} \|E_N(u_g)\|_2 = 0$.

\item For every free pmp action $\alpha : \Lambda \curvearrowright (Y, \nu)$, the co-induced action $\coInd(a, b)_\Lambda^\Gamma(\alpha)$ is mixing.
\end{itemize}
\end{theo}
Let $\rho = \coInd(a, b)_\Lambda^\Lambda(\alpha)$, $\sigma = \coInd(a, b)_\Lambda^\Gamma(\alpha)$ and assume that $\sigma$ is ergodic. The following properties hold true (see \cite{epstein}).
\begin{itemize}
\item [$(\ast)$] For any quotient map $q : Y \to Z$ from $\alpha : \Lambda \curvearrowright Y$ to a free pmp action $\Lambda \curvearrowright Z$, we have that 
$$\left \{ (x, (y_n)_{n \in \N}) : q \circ p_\rho(g^\sigma \cdot(x, (y_n)_{n \in \N})) = q \circ p_\rho((x, (y_n)_{n \in \N})) \right\}$$
is a $\mu \times \nu^\N$-null measurable subset, for all $g \in \Gamma \setminus \{1_\Gamma\}$.
\item [$(\ast \ast)$] For any $\rho(\Lambda)$-invariant Borel subset $U \subset X \times Y^\N$ of $\mu \times \nu^\N$-positive measure, the Borel map $p_\rho |U : U \to Y$ witnesses that $\alpha$ is a quotient of $\rho | U$.
\end{itemize}

Gaboriau and Lyons proved that given any nonamenable group $\Gamma$, there exist free pmp actions $a : \F_2 \curvearrowright (X, \mu)$ and $b : \Gamma \curvearrowright (X, \mu)$ such that $a$ is ergodic, $b$ is mixing and $\mathcal R(a, \F_2) \subset \mathcal R(b, \Gamma)$ (see Theorem \ref{measurable}). Epstein, Ioana, Kechris and Tsankov proved \cite[Theorem 3.11]{IKT} that the inclusion $\mathcal R(a, \F_2) \subset \mathcal R(b, \Gamma)$ can be chosen to satisfy the assumptions of Theorem $\ref{mixing}$.

\begin{theo}\label{measurable-bis}
Let $\Gamma$ be any nonamenable group. Then there exist free pmp actions $a : \F_2 \curvearrowright (X, \mu)$ and $b : \Gamma \curvearrowright (X, \mu)$ such that $a$ is ergodic, $b$ is mixing, $\mathcal R(a, \F_2) \subset \mathcal R(b, \Gamma)$ and $\lim_{g \to \infty} \|E_{L^\infty(X) \rtimes \F_2}(u_g)\|_2 =~0$.
\end{theo}

\section{Uncountably many non-OE actions}\label{uncountable}
\subsection{Separability vs.\ relative property (T)}

Recall that for an inclusion $\Lambda < \Gamma$ of countable discrete groups, the pair $(\Gamma, \Lambda)$ has the relative property (T) if for all $\varepsilon > 0$, there exist $\delta > 0$ and a finite subset $F \subset \Gamma$ such that if $\pi : \Gamma \to \mathcal U(\mathcal H)$ is a unitary representation and $\xi \in \mathcal H$ is a unit vector which satisfies $\|\pi(g)(\xi) - \xi\| < \delta$, for all $g \in F$, then there exists a $\pi(\Lambda)$-invariant vector $\eta \in \mathcal H$ such that $\|\eta - \xi\| < \varepsilon$. The pair $(\Z^2 \rtimes \SL_2(\Z), \Z^2)$ has the relative property (T) \cite{{kazhdan}, {margulis}}. More generally, for any nonamenable subgroup $\Gamma < \SL_2(\Z)$, the pair $(\Z^2 \rtimes \Gamma, \Z^2)$ has the relative property (T) \cite{burger}. 

Consider the action $\SL_2(\Z) \curvearrowright (\mathbf T^2, \lambda^2)$ defined by
$$g \cdot (z_1, z_2) = (g^{-1})^t \begin{pmatrix}
z_1 \\
z_2
\end{pmatrix}, \forall g \in \SL_2(\Z).$$
One checks that it is a free weakly mixing pmp action.
Realize $\F_2 < \SL_2(\Z)$ as a finite index subgroup, so that the pair $(\Z^2 \rtimes \F_2, \Z^2)$ has the relative property (T). Write $\alpha : \F_2 \curvearrowright (\mathbf T^2, \lambda^2)$ for the restriction.

The following result is due to Ioana \cite[Theorem 1.3]{ioana07}. It relies on a separability vs. (relative) property (T) argument, an idea that goes back to Connes \cite{connes-fundamental} and successfully used later on by Popa \cite{popa2001} and Gaboriau and Popa in \cite{gaboriau-popa}.

\begin{theo}\label{separability}
Let $\Gamma$ be any nonamenable group. Let $\mathcal F(\Gamma)$ be the class of free ergodic pmp actions $\sigma : \Gamma \curvearrowright (X, \mu)$ such that there exists a free pmp action $\rho : \F_2 \curvearrowright (X, \mu)$ for which the following hold:
\begin{enumerate}
\item $\mathcal R(\rho, \F_2) \subset \mathcal R(\sigma, \Gamma)$.
\item The action $\alpha : \F_2 \curvearrowright \mathbf T^2$ is a quotient of the action $\rho : \F_2 \curvearrowright X$ with quotient map $p_\rho : X \to \mathbf T^2$.
\item For all $g \in \Gamma \setminus \{1_\Gamma\}$, the Borel set $\{x \in X : p_\rho(\sigma(g)(x)) = p_\rho(x)\}$ is null.
\end{enumerate}
Let $\{\sigma_i : i \in \mathcal I\} \subset \mathcal F(\Gamma)$ be an uncountable set of mutually orbit equivalent actions. Then there exist an uncountable set $\mathcal J \subset \mathcal I$ and $\rho_j$-invariant measurable subsets $X_j \subset X$ of positive measure such that the actions $\{ \rho_j | X_j : j \in \mathcal J\}$ are mutually conjugate.
\end{theo}

\begin{proof}
By assumption, denote by $\mathcal R$ the unique pmp equivalence relation on $(X, \mu)$ (up to orbit equivalence) such that $\mathcal R = \mathcal R(\sigma_i, \Gamma)$, for all $i \in \mathcal I$. Note that for all $i \in \mathcal I$, $\mathcal R(\rho_i, \F_2) \subset \mathcal R$. Following \cite{feldman-moore}, define a Borel measure $\nu$ on $\mathcal R$ by
$$\nu(\mathcal W) = \int_X |\{y : (x, y) \in \mathcal W\}| {\rm d}\mu(x),$$
for every Borel subset $\mathcal W \subset \mathcal R$.

For all $i \in \mathcal I$, denote by $p_i : X \to \mathbf T^2$ the quotient map which witnesses that $\alpha : \F_2 \curvearrowright \mathbf T^2$ is a quotient of $\rho_i : \F_2 \curvearrowright X$. Regarding $a \in \Z^2$ as a character of $\mathbf T^2$, define $f_{a, i} = a \circ p_i \in L^\infty(X)$. One checks that for all $(a, g) \in \Z^2 \rtimes \F_2$ and $i \in \mathcal I$, $f_{g(a), i} = f_{a, i} \circ \rho_i(g^{-1})$. Then for all $i, j \in \mathcal I$, the map $\pi_{i, j} : \Z^2 \rtimes \F_2 \to \mathcal U(L^2(\mathcal R, \nu))$ defined by $\pi_{i, j}(a, g)(\xi)(x, y) = f_{a, i}(x) \overline{f_{a, j}(y)} \xi(\rho_i(g^{-1})(x), \rho_j(g^{-1})(y))$,
for all $(a, g) \in \Z^2 \rtimes \F_2$, $\xi \in L^2(\mathcal R, \nu)$, $(x, y) \in \mathcal R$, is a unitary representation.

Denote by $\Delta = \{(x, x) : x \in X\} \subset \mathcal R$ the diagonal. Note that $\mathbf 1_\Delta \in L^2(\mathcal R, \nu)$ and $\|\mathbf 1_\Delta\|_2 = 1$. One checks that for all $(a, g) \in \Z^2 \rtimes \F_2$, $i, j \in \mathcal I$,
$$
\|\pi_{i, j}(a, g)(\mathbf 1_\Delta) - \mathbf 1_\Delta\|_2^2 \leq  2 \|\mathbf 1_{\graph(\rho_i(g^{-1}))} - \mathbf 1_{\graph(\rho_j(g^{-1}))}\|_2 + 2 \|f_{a, i} \mathbf 1_\Delta - f_{a, j} \mathbf 1_\Delta \|_2.
$$
Since the pair $(\Z^2 \rtimes \F_2, \Z^2)$ has the relative property (T), with $\varepsilon = 1/2$, there exist $\delta > 0$, finite subsets $A \subset \Z^2$, $F \subset \F_2$ such that if $\pi : \Z^2 \rtimes \F_2 \to \mathcal U(\mathcal H)$ is a unitary representation and $\xi \in \mathcal H$ is a unit vector which satisfies $\|\pi(a, g)(\xi) - \xi\| < \delta$, for all $a \in A$ and $g \in F$, then there exists a $\pi(\Z^2)$-invariant vector $\eta \in \mathcal H$ such that $\|\eta - \xi\| < \varepsilon$. Since $\mathcal I$ is uncountable and $L^2(\mathcal R, \nu)$ is $\|\cdot\|_2$-separable, there exists an uncountable subset $\mathcal J \subset \mathcal I$, such that for all $i, j \in \mathcal J$,
\begin{eqnarray*}
\|f_{a, i} \mathbf 1_\Delta - f_{a, j} \mathbf 1_\Delta \|_2 & < & \delta^2/4, \forall a \in A \\
\|\mathbf 1_{\graph(\rho_i(g^{-1}))} - \mathbf 1_{\graph(\rho_j(g^{-1}))}\|_2 & < & \delta^2/4, \forall g \in F.
\end{eqnarray*}
Fix now $i, j \in \mathcal J$. Since $\|\pi_{i, j}(a, g)(\mathbf 1_\Delta) - \mathbf 1_\Delta\|_2 < \delta$, for all $(a, g) \in A \times F$, the relative property (T) gives a $\pi_{i, j}(\Z^2)$-invariant vector $\eta \in L^2(\mathcal R, \nu)$ such that $\|\eta - \mathbf 1_\Delta\|_2 \leq 1/2$. Hence, $\nu$-a.s. $\eta(x, y) = f_{a, i}(x) \overline{f_{a, j}(y)} \eta(x, y)$, for all $a \in \Z^2$. Since $\eta \neq 0$, the measurable subset $\mathcal W = \{(x, y) \in \mathcal R: f_{a, i}(x) = f_{a, j}(y), \forall a \in \Z^2\}$ satisfies $\nu(\mathcal W) > 0$. Next we claim that for $\mu$-almost every $x \in X$, there exists at most one $y \in X$ such that $(x, y) \in \mathcal W$. Assume this is not the case. Since $\mathcal R = \mathcal R(\sigma_j, \Gamma)$, one can find a measurable subset $Y \subset X$ of $\mu$-positive measure and $s \neq t \in \Gamma$, such that $(x, \sigma_j(s)(x)) \mbox{ and } (x, \sigma_j(t)(x)) \in \mathcal W$, for all $x \in Y$. In particular, we get $a\left(p_j(\sigma_j(s)(x))\right) = a \left( p_j(\sigma_j(t)(x)) \right)$, for all $a \in \mathbf Z^2, x \in Y$. Since characters separate points, it follows that $p_j(\sigma_j(s)(x)) = p_j(\sigma_j(t)(x))$, for all $x \in Y$. This clearly contradicts item $(3)$ in the statement of the Theorem. 

Define the measurable subset $X_i = \{x \in X : \exists! y \in X, (x, y) \in \mathcal W\}$. Since $\nu(\mathcal W) > 0$, the above claim yields $\mu(X_i) > 0$. If $(x, y) \in \mathcal W$, then $f_{a, i}(x) = f_{a, j}(y)$, for all $a \in \Z^2$ and hence $f_{g(a), i}(x) = f_{g(a), j}(y)$, for all $a \in \Z^2$, $g \in \F_2$.  Since $f_{g(a), i} = f_{a, i} \circ \rho_i(g^{-1})$, we get
\begin{equation}\label{relation}
(\rho_i(g)(x), \rho_j(g)(y)) \in \mathcal W, \forall g \in \F_2, \forall (x, y) \in \mathcal W.
\end{equation}
In particular, $X_i$ is a $\rho_i(\F_2)$-invariant measurable subset. Likewise, define $X_j = \{y \in X : \exists x \in X_i, (x, y) \in \mathcal W\}$. Then $X_j$ is a $\rho_j(\F_2)$-invariant measurable subset. Define $\phi : X_i \to X_j$ by $y = \phi(x)$ if and only if $(x, y) \in \mathcal W$. One checks that $\phi$ is a pmp Borel isomorphism. Finally, $(\ref{relation})$ shows that $\phi$ is a conjugacy between $\rho_i |X_i$ and $\rho_j | X_j$, that is, $\phi(\rho_i(g)(x)) = \rho_j(g)(\phi(x))$, for all $x \in X_i$, $g \in \F_2$.
\end{proof}

\subsection{A continuum of actions}
Let $\Gamma$ be any nonamenable group. Choose $a : \F_2 \curvearrowright (X, \mu)$ and $b : \Gamma \curvearrowright (X, \mu)$ according to Theorem \ref{measurable-bis}. Let $\pi : \F_2 \to \mathcal U(\mathcal H_\pi)$ be a unitary representation. Denote by $\gamma_\pi : \F_2 \curvearrowright (Z_\pi, \eta_\pi)$ the corresponding pmp Gaussian action (see \cite[Appendix E]{kechris} for more details).
\begin{itemize}
\item If $\pi_1$ and $\pi_2$ are unitarily equivalent, then $\gamma_{\pi_1}$ and $\gamma_{\pi_2}$ are conjugate.
\item If we denote by $\kappa(\gamma_\pi) : \F_2 \to \mathcal U(L^2(Z_\pi, \eta_\pi) \ominus \C1)$ the associated Koopman representation, we have $\pi \subset \kappa(\gamma_\pi)$.
\end{itemize}
Let $\alpha_\pi = \alpha \times \gamma_\pi : \F_2 \curvearrowright (\mathbf T^2 \times Z_\pi, \lambda^2 \times \eta_\pi)$ be the diagonal action. Observe that $\alpha_\pi$ is a free pmp action and $\alpha$ is a quotient of $\alpha_\pi$ via the quotient map $(y, z) \mapsto y$. Define the actions $\sigma_\pi := \coInd(a, b)_{\mathbf F_2}^\Gamma(\alpha_\pi)$ and $\rho_\pi := \coInd(a, b)_{\F_2}^{\F_2}(\alpha_\pi)$. Recall from Section \ref{induction} that $\sigma_\pi$ is mixing (see Theorem \ref{mixing}) and the following hold true:
\begin{enumerate}
\item $\mathcal R(\rho_\pi, \F_2) \subset \mathcal R(\sigma_\pi, \Gamma)$.
\item $\alpha$ is a quotient of $\rho_\pi$ with quotient map 
$$p_\pi :  X \times (\mathbf T^2 \times Z_\pi)^\N \ni (x, (y_n, z_n)_{n \in \N}) \mapsto y_0 \in \mathbf T^2.$$
\item For all $g \in \Gamma \setminus \{1_\Gamma\}$, the Borel set 
$$\left\{ (x, (y_n, z_n)_{n \in \N}) : p_\pi(g^{\sigma_\pi} \cdot (x, (y_n, z_n)_{n \in \N})) = p_\pi((x, (y_n, z_n)_{n \in \N})) \right\}$$
is $\mu \times (\lambda^2 \times \eta_\pi)^\N$-null (by Condition $(\ast)$ from Section \ref{induction}).
\end{enumerate}
The last result of this text is \cite[Theorem 5]{IKT}. We point out that it was first obtained by Ioana \cite[Section 3]{ioana07} when $\F_2 < \Gamma$ and then extended by Epstein \cite{epstein} when $\F_2 <_{\ME} \Gamma$ but without the mixing property.

\begin{theo}\label{nonOE}
Let $\Gamma$ be any nonamenable group. Then $\Gamma$ admits uncountably many non-orbit equivalent free mixing pmp actions.
\end{theo}

\begin{proof}
Let $\mathcal I_0$ be an uncountable set of pairwise non-isomorphic irreducible representations of $\F_2$ (see \cite{szwarc}). Denote by $(\mathcal U, \tau)$ the standard Borel probability space $(X \times (\mathbf T^2 \times Z)^\N, \mu \times (\lambda^2 \times \eta)^\N)$.  By contradiction, assume that there exists an uncountable subset $\{\sigma_\pi : \pi \in \mathcal I\} \subset \mathcal F(\Gamma)$ of mutually orbit equivalent actions. By Theorem \ref{separability}, there exists an uncountable subset $\mathcal J \subset \mathcal I$ and $\rho_\pi$-invariant Borel subsets $\mathcal U_\pi \subset \mathcal U$ of $\tau$-positive measure such that the actions $\{ \rho_\pi |\mathcal U_\pi : \pi \in \mathcal J\}$ are mutually conjugate. By Condition $(\ast \ast)$ from Section \ref{induction}, we know that $\alpha \times \gamma_\pi$ is a quotient of $\rho_\pi | \mathcal U_\pi$. Fix now $\pi_0 \in \mathcal J$. For all $\pi \in \mathcal J$, we have
$$\pi \subset \kappa(\gamma_\pi) \subset \kappa(\alpha \times \gamma_\pi) \subset \kappa(\rho_\pi |\mathcal U_\pi) \cong \kappa(\rho_{\pi_0}| \mathcal U_{\pi_0}) \subset \kappa({\rho_{\pi_0}}).$$
Then the separable unitary representation $\kappa({\rho_{\pi_0}})$ contains uncountably many pairwise non-isomorphic irreducible subrepresentations $\pi \in \mathcal J$, which is a contradiction.
\end{proof}

\bibliographystyle{plain}
\end{document}